\documentclass[11pt,english,reqno]{smfart}
\usepackage[letterpaper,left=1in,right=1in,top=1.5in,bottom=1.5in]{geometry}

\usepackage[T1]{fontenc}
\usepackage[english,francais]{babel}

\usepackage{amssymb,url,xspace,smfthm}

\newcommand{\SMF}{Soci\'et\'e Ma\-th\'e\-Ma\-ti\-que de France}
\newcommand{\BibTeX}{{\scshape Bib}\kern-.08em\TeX}
\newcommand{\T}{\S\kern .15em\relax }
\newcommand{\AMS}{$\mathcal{A}$\kern-.1667em\lower.5ex\hbox
        {$\mathcal{M}$}\kern-.125em$\mathcal{S}$}

\usepackage[hidelinks]{hyperref}
\usepackage[
]{cleveref}

\usepackage{mathtools}

\title[On $q$-series and local factors]{On $q$-series and the moment problem\\ associated to  local factors}
\date {}
\author{Alain Connes}
\address{Coll\`ege de France\\3 Rue d'Ulm\\75005 Paris, France}
\address{I.H.E.S., France}
\email{\url{alain@connes.org}}
\author{Caterina Consani}
\address{Department of Mathematics\\Johns Hopkins University\\Baltimore, MD 21218, USA}
\email{\url{cconsan1@jhu.edu}}
\author{Henri Moscovici}
\address{Department of Mathematics\\Ohio State University\\Columbus, OH 43210, USA}
\email{\url{moscovici.1@osu.edu}}

\subjclass{\href{http://www.ams.org/mathscinet/msc/msc2020.html?t=11B65&btn=Current}{11B65},
\href{http://www.ams.org/mathscinet/msc/msc2020.html?t=33D15&btn=Current}{33D15},
\href{http://www.ams.org/mathscinet/msc/msc2020.html?t=33C45&btn=Current}{33C45},
\href{http://www.ams.org/mathscinet/msc/msc2020.html?t=42C05&btn=Current}{42C05},
\href{http://www.ams.org/mathscinet/msc/msc2020.html?t=47B36&btn=Current}{47B36}
\href{http://www.ams.org/mathscinet/msc/msc2020.html?t=05A10&btn=Current}{05A10}
\href{http://www.ams.org/mathscinet/msc/msc2020.html?t=05A19&btn=Current}{05A19}}

\keywords{Semilocal trace formula, q-series,  Orthogonal polynomials,  Jacobi matrix, Catalan numbers}

\usepackage{amsmath}
\usepackage{amscd}
\usepackage{amsbsy}
\usepackage{amssymb}
\usepackage{verbatim}
\usepackage{eufrak}
\usepackage{eucal}
\usepackage{microtype}
\usepackage{mathrsfs}
\usepackage{amsthm}
\usepackage{stmaryrd}
\usepackage[all,cmtip]{xy}
\usepackage{tikz}
\usetikzlibrary{matrix,arrows}
\usepackage[abbrev,lite,alphabetic]{amsrefs}
\usepackage{dsfont}

\usepackage{color}
\definecolor{todo}{rgb}{1,0,0}
\definecolor{conditional}{rgb}{0,1,0}
\definecolor{e-mail}{rgb}{0,.40,.80}
\definecolor{reference}{rgb}{.20,.60,.22}
\definecolor{mrnumber}{rgb}{.80,.40,0} 
\definecolor{citation}{rgb}{0,.40,.80} 

\newcommand{\ie}{{\it i.e.\/}\ }
\newcommand{\eg}{{\it e.g.\/}\ }
\newcommand{\qqq}{\,,\,~\forall}
\newcommand{\cf}{{\it cf.}}
\newcommand{\opcit}{{\it op.cit.\/}\ }

\newcommand{\GL}{{\rm GL}}

\def\End{{\rm End}}

\def\GL{{\rm GL}}

\def\id{{\rm id}}

\def\Tr{{\rm Tr}}
\def\tr{{\rm Tr}}

\def\A{{\mathbb A}}

\def\N{{\mathbb N}}

\def\Q{{\mathbb Q}}
\def\R{{\mathbb R}}
\def\Z{{\mathbb Z}}

\def\T{{\mathbb T}}

\def\Tr{{\rm Tr}}
\def\tr{{\rm tr}}

\def\cA{{\mathcal A}}
\def\cB{{\mathcal B}}
\def\cC{{\mathcal C}}
\def\cD{{\mathcal D}}
\def\cE{{\mathcal E}}

\def\cH{{\mathcal H}}

\def\cN{{\mathcal N}}

\def\cP{{\mathcal P}}

\def\cR{{\mathcal R}}
\def\cS{{\mathcal S}}

\def\qqq{\,,\,~\forall}

\def\ccyc{\Z^{\rm cyc}}

\def\id{{\mbox{Id}}}

\def\dim{{\mbox{dim}}}
 \def\scal2{{\mathscr S}}

\def\sr0{{\cS^{\rm ev}_0}}

\def\sar0{{\cS_0(\A_\Q)}}

\def\scal{{\mathbb S}}

\newcommand{\nil}[1]{}

\theoremstyle{plain}
\newtheorem{thm}{Theorem}[section]

\newtheorem{cor}[thm]{Corollary}
\newtheorem{lem}[thm]{Lemma}

\newtheorem{fact}[thm]{Fact}

\begin{document}
\def\smfbyname{}

\begin{abstract}
We investigate the moment problem and Jacobi matrix associated-- by  the operator theoretic framework of the semilocal trace formula-- to  
each finite set $S$ of places of $\Q$ containing the archimedean place.  The  measure is given  by the absolute value squared of the product over $S$ of local factors restricted to the critical line. 
We  treat  the case $S=\{p,\infty\}$, where a single prime $p$ is adjoined to the archimedean place. We find that all the key ingredients such as the moments, the orthogonal polynomials and the Jacobi matrices  can be expressed as power series in terms of the parameter $q:=1/p$. We show  that the series which appear for the moments themselves are Lambert series.  The study of the $q$-series for the coefficients of the Jacobi matrix, and for the associated orthogonal polynomials reveals an intriguing integrality result : all those coefficients belong to the ring $\Z[\frac{1}{\sqrt{2}}]$ obtained by adjoining 
$1 / \sqrt{2}$ to the ring of integers. The main result of this paper is the conceptual explanation of this integrality property using Catalan numbers.
\end{abstract}

\maketitle

\tableofcontents

\section{Introduction}
 The semilocal trace formula \cite{C98,CMbook} gives an operator theoretic framework suitable to handle Weil positivity as in \cite{CCweil}, as well as the spectral realization of the ultraviolet behavior of the zeros of the Riemann zeta function as in \cite{CM}. Using the notion of cyclic pair applied to the scaling operator,
we showed in \cite{c2m} that this gives rise, for 
each finite set $S$ of places of $\Q$ containing the archimedean place, to a determinate moment problem whose associated measure is governed by the product over $S$ of local factors restricted to the critical line. In the  purely archimedean case, \ie when $S=\{\infty\}$, the moments are given by rational numbers whose numerators form a classical sequence of integers (known as $A126156$) and with denominators given by powers of $2$.
The investigation of these moment problems in the general $S$-semilocal  case is meant to provide the extension of the Jacobi picture for the scaling operator-- treated  in \cite[\S 5]{c2m} at the archimedean place -- to the general  case.\newline 
In the present paper we  treat the case  where a single prime $p$ is adjoined to the archimedean place, \ie  the case $S=\{\infty, p\}$. We find  that all the essential ingredients for the Jacobi picture, 
such as the moments, the orthogonal polynomials and the Jacobi matrix, can be expressed in terms 
of power series in terms of the parameter $q=1/p$. The series which appear in the expression of the moments themselves are Lambert series.   
Furthermore, the calculation of the $q$-series arising as coefficients of the Jacobi matrix and of the associated orthogonal polynomials reveals an intriguing integrality result, namely 
the fact that all those coefficients belong to the ring $\Z[\frac{1}{\sqrt{2}}]$. The main outcome of the present paper is a conceptual explanation of this integrality property. \newline 
  The outline of the paper is as follows. In \S \ref{sectorthopoly} we recall a few basic 
 facts about orthogonal polynomials, including an explanation of their spectral meaning 
 (Proposition \ref{spec}) for which we did not find an adequate reference in the literature,
 although its proof is based on a classical (and often rediscovered \cf\ \cite{berg}) result of A.C. Aitken about the inverse of a Hankel matrix. 
 We then proceed in \S \ref{sectsemi} to the computation of the orthogonal polynomials associated 
 to the moment problem corresponding the general semilocal case. After showing that the moment problem remains determinate (Proposition \ref{momentdet}) we find the generating function of
the moments (Proposition \ref{momentssm}). 

For the rest of the paper we restrict ourselves to the case of a single added prime, \ie
$S=\{\infty, p\}$. The Lambert series make their appearance in \S \ref{lambert} in the
explicit expression of the moments (Proposition \ref{expth3}).  Concrete computer-aided
calculations for the deformation of the Jacobi matrix in the presence of an added prime reveal 
a surprising integrality pattern (Fact \ref{keyfact}) and the major part of the paper is devoted
to its ramifications and rigorous proof. Using Aitken's result quoted above and
the Catalan numbers, we prove in \S \ref{secinteg} a general integrality property (Theorem 5.9) 
which, besides proving Fact \ref{keyfact}, leads the way for obtaining in the subsequent section
\S \ref{eigenv} a companion integrality property (Theorem \ref{eigenthm}). The latter property
carries over to the coefficients of all $q$-series involved in the process of $q$-deforming 
the Jacobi matrix and of the orthogonal polynomials. It also holds for the deformation of the
diagonal parameters $(d_n)$ which single out the metaplectic representation in the Archimedean case (\cf \, \cite[\S 5]{c2m}), and in the final section \S \ref{diag} 
we compute the first terms of the $q$-expansion of the coefficients $d_n(q)$ and show that they are expressed in terms of hypergeometric functions $Y:=\, _2F_1\left(1,n+\frac{1}{2};n+1;-1\right)$.

\section{Moments and orthogonal polynomials}\label{sectorthopoly}
We review in this section basic properties of orthogonal polynomials, first their meaning as characteristic polynomials in \S \ref{spectralpol}, and then the special properties of the even case in \S \ref{evenpol}.
 \subsection{Spectral meaning of orthogonal polynomials}\label{spectralpol}
  Let $d\mu$ be a probability measure on $\R$. Orthogonal polynomials are defined by orthogonalization process from the powers $x^n$ of the variable and they form an orthonormal basis of the Hilbert space $L^2(\R, d\mu)$. There is however another understanding of their role, showing that their zeros compute the spectrum of the compression of the operator of multiplication by $x$ on the subspace of $L^2(\R, d\mu)$ given by polynomials of degree $< n$.	 \newline
  We  work in the general framework of orthogonal polynomials. We start from the formula (up to normalization) which gives the orthogonal polynomial $P_n(x)$ as a determinant  (see \cite{Szego}, (2.2.6) page 27). We let $c(n):=\int x^n d\mu(x)$ be the moments of the measure $d\mu$ and $D_n$ be the determinant of the Hankel matrix $\cD(n)$ with entries $c(i+j)$ for $i,j\in \{0,\ldots, n\}$. One then has, by (2.2.9) page 27, of \opcit that the orthogonal polynomial $P_n(x)$ is given by 
\begin{equation}\label{detorthopol}
	P_n(x)=(-1)^n(D_nD_{n-1})^{-1/2}\det(\cD'(n-1)-x\ \cD(n-1))
\end{equation}
where $\cD'(n)$ is the Hankel matrix with entries $c(i+j+1)$ for $i,j\in \{0,\ldots, n\}$.		More specifically, and with \eg   $n=4$, one starts from (2.2.6) page 27, of \opcit
   $$
  P_4(x)=(D_4D_{3})^{-1/2}\det \left(
\begin{array}{ccccc}
 c(0) & c(1) & c(2) & c(3) & c(4) \\
 c(1) & c(2) & c(3) & c(4) & c(5) \\
 c(2) & c(3) & c(4) & c(5) & c(6) \\
 c(3) & c(4) & c(5) & c(6) & c(7) \\
 1 & x & x^2 & x^3 & x^4 \\
\end{array}
\right)
$$
One then reduces the last line to $(1,0,0,\ldots)$ by replacing the $j$-th column $\cC_j$ by   $\cC_j-x\cC_{j-1}$ to get 
$$
 P_4(x)=(D_4D_{3})^{-1/2}\det\left(
\begin{array}{cccc}
 c(1)-c(0) x & c(2)-c(1) x & c(3)-c(2) x & c(4)-c(3) x \\
 c(2)-c(1) x & c(3)-c(2) x & c(4)-c(3) x & c(5)-c(4) x \\
 c(3)-c(2) x & c(4)-c(3) x & c(5)-c(4) x & c(6)-c(5) x \\
 c(4)-c(3) x & c(5)-c(4) x & c(6)-c(5) x & c(7)-c(6) x \\
\end{array}
\right)
$$ 
The next Proposition holds in the general case of orthogonal polynomials.
\begin{prop}\label{orthopol1} The orthogonal polynomial $P_n$ is, up to a multiplicative constant, the characteristic polynomial of the compression $\Pi_{n-1}X\ \Pi_{n-1}$ of the operator $X$ of multiplication by $x$ on the subspace $\cH_{n-1}$ of $L^2(\R, d\mu)$ given by polynomials of degree $< n$.	
\end{prop} \label{spec}
\proof Let $P_\ell(x)=\sum b_{k,\ell}\, x^k$ be the orthonormal polynomials.
By the result of A.C.~Aitken (see \cite{berg}), the inverse of the Hankel matrix $\cD(n)$ is equal to $\cB_n\, \cB_n^*$ where $\cB_n$ is  the upper triangular matrix with matrix entries $b_{k,\ell}$ and $\cB_n^*$ is its adjoint
$$
\cB=\left(
\begin{array}{ccccc}
 b_{0,0} & b_{ 0,1 } & b_{ 0,2 } & b_{ 0,3 } & b_{ 0,4 } \\
 0 & b_{ 1,1 } & b_{ 1,2 } & b_{ 1,3 } & b_{ 1,4 } \\
 0 & 0 & b_{ 2,2 } & b_{ 2,3 } & b_{ 2,4 } \\
 0 & 0 & 0 & b_{ 3,3 } & b_{ 3,4 } \\
 0 & 0 & 0 & 0 & b_{ 4,4 } \\
\end{array}
\right)
$$
$$
\cB^*=\left(
\begin{array}{ccccc}
 \overline b_{ 0,0 } & 0 & 0 & 0 & 0 \\
 \overline b_{ 0,1 } & \overline b_{ 1,1 } & 0 & 0 & 0 \\
 \overline b_{ 0,2 } & \overline b_{ 1,2 } & \overline b_{ 2,2 } & 0 & 0 \\
 \overline b_{ 0,3 } & \overline b_{ 1,3 } & \overline b_{ 2,3 } & \overline b_{ 3,3 } & 0 \\
 \overline b_{ 0,4 } & \overline b_{ 1,4 } & \overline b_{ 2,4 } & \overline b_{ 3,4 } & \overline b_{ 4,4 } \\
\end{array}
\right)
$$
In fact one has the identity
\begin{equation}\label{invbbstar}
	\cB_n^*\ \cD(n) \ \cB_n=\id
\end{equation}
since the entry $(i,\ell)$ of the  product matrix is given by 
$$
\sum_{j,k} \overline b_{ j,i }c(j+k)\ b_{ k,\ell }=\int \overline{P_i(x)}P_\ell(x)d\mu(x)=\delta_i^\ell
$$
 One  then gets 
$$
	\cB_n^*\ (\cD'(n)-z\ \cD(n)) \ \cB_n=\cB_n^*\ \cD'(n) \ \cB_n-z\ \id
		$$
and it remains to show that the matrix $\cB_n^*\ \cD'(n) \ \cB_n$ is the matrix of the compression of the operator of multiplication by $x$ on the subspace of $L^2(\R, d\mu)$ given by polynomials of degree $\leq n$. Its entry $(i,\ell)$  is given by 
$$
\sum \overline b_{ j,i }c(j+k+1)\ b_{ k,\ell }=\int \overline{P_i(x)}P_\ell(x)x\ d\mu(x)=\langle xP_i\vert P_\ell\rangle
$$
which gives the required result, since the $P_\ell$, for $0\leq \ell\leq n$, form an orthonormal basis of the subspace of $L^2(\R, d\mu)$ given by polynomials of degree $\leq n$.	\endproof 
\subsection{The even case}\label{evenpol}
 We assume that the measure $d\mu$ is even, \ie that $d\mu(-x)=d\mu(x)$. The  moments  $c(n):=\int x^n d\mu(x)$ are $0$ for odd $n$. The orthogonal polynomials $P_\ell$ for $0\leq \ell\leq n$, form an orthonormal basis of the subspace $\cH_n$ of $L^2(\R, d\mu)$ given by polynomials of degree $\leq n$, uniquely specified if one requires that the leading coefficient of $P_\ell$ is positive. The operator $X$ of multiplication by the variable $x\in \R$ is then given by a Jacobi matrix and one has non-zero scalars $a_n:=\langle X\ P_n\vert P_{n+1}\rangle$ such that 
 \begin{equation}\label{dxij}
X\ P_n=a_{n-1}\ P_{n-1}+a_{n}\ P_{n+1}, \ \ a_n\neq 0, 
\end{equation}
so that the  matrix of $X$ has the following form
$$
X=\left(
\begin{array}{cccccc}
 0 & a_0 & 0 & 0 & 0 & \ldots \\
  a_0 & 0 & a_1 & 0 & 0 & \ldots \\
 0 &   a_1 & 0 & a_2 & 0 & \ldots \\
 0 & 0 &  a_2 & 0 & a_3 & \ldots \\
 0 & 0 & 0 &   a_3 & 0 & \ldots\\
 \ldots & \ldots & \ldots & \ldots & \ldots & \ldots \\
\end{array}
\right)
$$
Since the odd moments vanish, the Hankel matrix $\cD(n)$ has the following form
$$
\cD(4)=\left(
\begin{array}{ccccc}
 c(0) & 0 & c(2) & 0 & c(4) \\
 0 & c(2) & 0 & c(4) & 0 \\
 c(2) & 0 & c(4) & 0 & c(6) \\
 0 & c(4) & 0 & c(6) & 0 \\
 c(4) & 0 & c(6) & 0 & c(8) \\
\end{array}
\right)
$$
and can be decomposed as the direct sum of the restrictions to   the two subspaces $\cH_{n}^\pm\subset \cH_n$ of vectors whose all odd coordinates (resp. even) vanish.
$$
\cD(4)=\left(
\begin{array}{ccc}
 c(0) & c(2) & c(4) \\
 c(2) & c(4) & c(6) \\
 c(4) & c(6) & c(8) \\
\end{array}
\right)\bigoplus \left(
\begin{array}{cc}
 c(2) & c(4) \\
 c(4) & c(6) \\
\end{array}
\right)
$$
The classical formula for the coefficients $a_n^2$ in terms of $D_n:=\det(\cD_n)$ is (see \cite{Szego})
 \begin{equation}\label{an2}
 a_n^2=\frac{D_{n-1}\ D_{n+1}}{D_n^2}
  \end{equation}
  In our case the odd moments vanish, and one can then replace \eqref{an2} by the simpler formula
  \begin{equation}\label{an2bis}
 a_n^2=\frac{N_{n+1}}{D_n}
  \end{equation}
  where the term $N_{n+1}$ is the determinant of the matrix $\cN_{n+1}$, defining a minor,  obtained by deleting the penultimate line and column of the matrix $\cD_{n+1}$. In order to prove \eqref{an2bis},
   we shall decompose  the matrices $\cD_n$ and $\cN_{n}$ using the decomposition  $\cH_{n}=\cH^+_{n}\oplus \cH^-_{n}$. \begin{prop} \label{evenmoments} Assume that all odd moments vanish.  \newline
$(i)$~Let $\cD_n^\pm$ be the restriction of the Hankel matrix to $\cH_{n}^\pm$. One has 
\begin{equation}\label{an2bis}
D_n=\det(\cD_n)= \det(\cD_n^+)\det(\cD_n^-)
  \end{equation}
  $(ii)$~Let $N_{n}$ is the determinant of the matrix $\cN_{n}$. One has 
	\begin{equation}\label{nndet}
N_n=\det(\cN_n)= \begin{cases} \det(\cD_n^+)\det(\cD_{n-2}^-)\ \textit{if} \ n \ \textit{even} \\
	\det(\cD_{n-2}^+)\det(\cD_{n}^-)\ \textit{if} \ n \ \textit{odd}
\end{cases}
  \end{equation}
  $(iii)$~The coefficients $a(n)^2$ are given by 
  \begin{equation}\label{nndet1}
a_n^2= \begin{cases} \left(\det(\cD_{n+1}^+)/\det(\cD_{n}^+)\right)\left(\det(\cD_{n-1}^-)/\det(\cD_{n}^-)\right)\ \textit{if} \ n \ \textit{odd} \\
	\left(\det(\cD_{n-1}^+)/\det(\cD_{n}^+)\right)\left(\det(\cD_{n+1}^-)/\det(\cD_{n}^-)\right)\ \textit{if} \ n \ \textit{even}
\end{cases}
  \end{equation}
  $(iv)$~Let $C_n:=\det(\cD_n^+)$ if $n$ is even and $C_n:=\det(\cD_n^-)$ if $n$ is odd. One then  has for all $n$
   \begin{equation}\label{nndet2}
   	a_n^2= \left(C_{n-2}/C_n\right)\times \left(C_{n+1}/C_{n-1}\right)
   \end{equation}
   $(v)$~One has \eqref{an2bis} for all $n$.  
  \end{prop}
\proof $(i)$~Follows from the decomposition $\cD_{n}=\cD_{n}^+\oplus \cD_{n}^-$.\newline
$(ii)$~Deleting the penultimate line and column means that one restricts the quadratic form to the subspace defined by the condition $\xi_{n-1}=0$. This restriction does not affect the space $\cH_{n}^\epsilon$ where $\epsilon= (-1)^{n}$ since all vectors in that subspace fulfill $\xi_{n-1}=0$. It replaces $\cH_{n}^{-\epsilon}$ by $\cH_{n-2}^{-\epsilon}$ since it nullifies the top component of the vectors. In fact one even has $\cH_{n-2}^{-\epsilon}=\cH_{n-3}^{-\epsilon}$.\newline
$(iii)$~Assume first that $n$ is odd. One then has  $\cH_{n}^+= \cH_{n-1}^+$ and $\cH_{n+1}^-= \cH_{n}^-$, thus, using $(i)$, 
$$
\det(\cD_{n+1})/\det(\cD_{n})=\det(\cD_{n+1}^+)/\det(\cD_{n}^+)
$$
and 
$$
\det(\cD_{n-1})/\det(\cD_{n})=\det(\cD_{n-1}^-)/\det(\cD_{n}^-)
$$
which gives the required result. The proof when $n$ is even is similar.\newline
$(iv)$~This follows from $(iii)$, \ie \eqref{nndet1}, since for $n$ even one has 
$\det(\cD_{n-1}^+)=\det(\cD_{n-2}^+)=C_{n-2}$ and $\det(\cD_{n}^-)=\det(\cD_{n-1}^-)=C_{n-1}$, while for $n$ odd one has $\det(\cD_n^+)=\det(\cD_{n-1}^+)=C_{n-1}$ and $\det(\cD_{n-1}^-)=\det(\cD_{n-2}^-)=C_{n-2}$. \newline
$(v)$~Follows from $(ii)$ and $(iii)$.
\endproof
By construction the $a_n^2$ are rational functions of the moments of the form
\begin{equation}\label{ansq}
 a_n^2=R_n((c(2j)_{0\leq j\leq n+1})
\end{equation}
The first instances are 
$$
R_1=\frac{c(4)}{c(2)}-\frac{c(2)}{c(0)},\ \ R_2=\frac{c(0) \left(c(4)^2-c(2) c(6)\right)}{c(2)^3-c(0) c(2) c(4)},
$$
$$
R_3=\frac{c(2) \left(c(4)^3-(2 c(2) c(6)+c(0) c(8)) c(4)+c(0) c(6)^2+c(2)^2 c(8)\right)}{\left(c(2)^2-c(0) c(4)\right) \left(c(2) c(6)-c(4)^2\right)}
$$
$$
R_4=-\frac{\left(c(2)^2-c(0) c(4)\right) \left(-c(6)^3+(2 c(4) c(8)+c(2) c(10)) c(6)-c(2) c(8)^2-c(4)^2 c(10)\right)}{\left(c(4)^2-c(2) c(6)\right) \left(c(4)^3-(2 c(2) c(6)+c(0) c(8)) c(4)+c(0) c(6)^2+c(2)^2 c(8)\right)}
$$
These expressions are homogeneous of degree $0$.
 Thus one does not have to normalize the moments in the computations.

\section{Moments in the semilocal case }\label{sectsemi}
\subsection{Moments in the Archimedean case}
We begin by recalling \cite[Proposition 5.2]{c2m} according to which 
in the Archimedean case   
 the normalized moments $c_n$ of the measure $d\mu =(2 \pi)^{-3/2}\vert\Gamma(\frac 14+\frac 12 is)\vert^2ds$ are determined by the identity 
\begin{equation}\label{moments1}
\sum c_n\frac{(ix)^n}{n!}=\sqrt{\frac{2}{e^{x}+e^{-x}}} .
\end{equation}	
In particular $c_n=0$ for $n$ odd, and for $n$ even the first the first ones are 
\begin{equation}\label{unmoments}
c_0=1,\ c_2= \frac{1}{2}, \ c_4= \frac{7}{4},\ c_6= \frac{139}{8},\ c_8= \frac{5473}{16},\ c_{10}= \frac{357721}{32}
\end{equation}
What one finds in general is that the denominators of the moments are $2^n$ for $c_{2n}$  so that they are determined by a sequence of integers. The first ones are 
$$
1,\ 7,\ 139,\ 5473,\ 357721,\ 34988647,\ 4784061619,\  871335013633,\ 203906055033841,$$ $$59618325600871687,\ 21297483077038703899,\ 9127322584507530151393, \ldots $$
We define the un-normalised (even) moments $c_0$ by $$c_0(2k):=(-1)^k\left(\partial_x^{2k}\frac{1}{\sqrt{e^{-x}+e^x}}\right)_{x=0}$$
The first values are
$$
\left\{\frac{1}{\sqrt{2}},\frac{1}{2 \sqrt{2}},\frac{7}{4 \sqrt{2}},\frac{139}{8 \sqrt{2}},\frac{5473}{16 \sqrt{2}},\frac{357721}{32 \sqrt{2}}\right\}.$$

\subsection{Determinacy of moment problem} 
Passing to the general semilocal case, we let now $S$ be a finite set of places with $\infty\in S$,
and consider the measure $\ d\mu_S(s):=\Big\vert\prod_{v\in S} L_p (\frac 12-is)\Big\vert^2\, ds$ \ 
where \, $  L_\infty(\frac 12-is)  =\pi^{-1/4}  \Gamma(\frac 14 + \frac 12 is) $.

\begin{prop} \label{momentdet}
 The moment problem for the sequence 
 $\displaystyle \big\{c_n^S:=\int s^n d\mu_S(s)\big\}_{n \geq 0}$ 
 is determinate. 
\end{prop} 
\proof By Carleman's Criterion (see \cite[Theorem 4.3]{sm}), a Hamburger problem 
with moments  $\gamma_n:=\int_\R t^n d\mu(t)$ is determinate if \
 $\displaystyle \sum_{n \geq 1} \gamma_{2n}^{-\frac{1}{2n}} = \infty$.
 We shall actually show that 
\begin{align} \label{carleman2}
\exists \,  M>0 \quad \text{such that}  \quad \gamma_{2n} \leq M^n (2n)! \, ,\quad \forall n \in \N ,
\end{align} 
which implies Carleman's condition (see \cite[Corollary 4.10]{sm}).  \newline

Since excluding the Archimedean place the product
$\, \Big\vert\prod_{v\in S \setminus\{\infty\}} L_v(\frac 12-is)\Big\vert^2$ is bounded, it suffices to
verify \eqref{carleman2} for the sequence of Archimedean moments $\{c_n\}_{n \geq 0}$. 
Furthermore,  the ratio 
$\displaystyle \vert \Gamma(\frac 14 +  is)\vert ^2 / \vert \Gamma(\frac 12 +  is)\vert ^2$ is 
bounded too, since 
$$
\vert \Gamma(\frac 14 +  is)\vert / \vert \Gamma(\frac 12 +  is)\vert \sim |s|^{-\frac 14} \quad 
\text{as} \quad |s| \to \infty 
$$
(see https://dlmf.nist.gov/5.11 [(iii) Ratios)]. Thus, it remains to check that  
\begin{align} \label{carleman3}
\int_0^\infty s^{n}  \vert \Gamma(\frac 12 +  i\frac s2)\vert ^2  ds \leq M^n n! \, ,\quad \forall n \in \N 
\end{align} 
for some $M>0$. \,
Using the formula \, $\displaystyle  \vert \Gamma(\frac 12 +  i t)\vert ^2 = \frac{\pi}{\cosh \pi t}$, \, one has
\begin{align*}
\int_0^\infty s^{n}  \vert \Gamma(\frac 12 +  i\frac s 2)\vert ^2  ds& = 
2^{n+1}\int_0^\infty t^n \vert \Gamma(\frac 12 +  i t)\vert ^2  dt =
2^{n+1}\int_0^\infty t^{n} \frac{2\pi}{e^{\pi t} + e^{-\pi t}} dt \\
&\leq  2^{n+2}\pi \int_0^\infty t^{n}  e^{-\pi t} dt = 2^{n+2}\pi^{-n} n! ,
\end{align*}
which achieves the proof.

\subsection{Generating function of moments in the semilocal case}
The notation used in this section is the same as in \cite{c2m}.
Let $S$ be a finite set of places, $\infty\in S$.   
 and $\A_{S}$ 
the locally compact ring
$
 \A_{S}=\prod_{v\in S} \Q_v.
$
This ring contains $\Q$ as a subring using the diagonal embedding. Let
$\Q_S$ denote the subring of $\Q$ given by rational numbers whose
denominator only involves primes $p \in S$. In other words,
$
 \Q_S=\{q\in \Q\,|\, |q|_v\leq 1\,,\ \forall v\notin S\}\,.
$
The group $\Gamma:=\Q^*_S$ of invertible elements of the ring $\Q_S$ is of
the form
\begin{equation}\label{GL1QS}
\Gamma=\GL_1(\Q_S)= \{ \pm p_1^{n_1} \cdots p_k^{n_k} \, :\,  p_j
\in S \setminus\{ \infty \} \,,\, n_j\in \Z\}.
\end{equation}
 The semilocal ad\`ele  class space  $X_{S}$ is by definition the quotient
\begin{equation}\label{XQS}
X_{S}:=\A_{S}/\Gamma,
\end{equation}
and  we let $\pi_S:\A_{S}\to X_S$ be the canonical projection. Let $R_S$ be the maximal compact subring of $\prod_{p\neq \infty}\Q_p$. For $f\in L^2(\R)^{ev}$ let $\eta_S(f)$ be the class of the function $1_{R_S}\otimes f$ in $L^2(X_S)$. More precisely one has, ignoring the summation over $\pm 1$ since $f$ is even, and letting $\Gamma_+:=\Gamma\cap \Q_+$,
 \begin{equation}\label{etas}
 \eta_S(f)(x):=\sum _{\Gamma_+}\,f(\gamma \,\tilde x)\qqq \tilde x \mid \pi_S(\tilde x)=x.
 \end{equation}
  The action $\exp(it\scal)$ of the scaling group on $L^2(\R)^{ev}$ fulfills, for $f\in L^2(\R)^{ev}$ and $t\in \R$,
 \begin{equation}\label{scalrule}
  \exp(it\scal)\eta_S(f)=\eta_S\left(\exp(it\scal)f\right), \ \ \exp(it\scal)f(x)=\lambda^{1/2}f(\lambda x), \ \lambda=e^t
  \end{equation}
We let $\xi_\infty=h_0\in L^2(\R)^{ev}$ be the vector of norm $1$ given by 
$$
h_0(x):=2^{\frac 14}\,\exp(-\pi x^2)\qqq x\in \R.
$$
\begin{prop} \label{momentssm}
Let $S$ be a finite set of places, $\infty\in S$, $\xi_S=\eta_S(\xi_\infty)$, one has 
\begin{equation}\label{expth}
\langle\exp(it\scal)\xi_S\vert \xi_S\rangle =\sum_{a/b\in \Gamma_+} \frac{\sqrt 2}{\sqrt{a^2e^t+b^2e^{-t}}}
\end{equation}
where the sum is over relatively prime pairs $(a,b)$ of elements of $\Gamma_+\cap\Z$.	
\end{prop}
\proof Let $R=R_S$, $\xi=1_R\otimes \xi_\infty$.  The inner product in $L^2(X_S)$ gives 
$$
\langle\exp(it\scal)\xi\vert \xi\rangle=\sum_{\Gamma_+}\langle \theta(\gamma)\xi\vert \exp(-it\scal)\xi \rangle
$$ 
With $\gamma=\frac ab= \prod p_v^{n_v} $, one has 
$$
\theta(\gamma)(1_R\otimes \xi_\infty)(y,x)=2^{\frac 14}\,\prod 1_{\Z_v}(p_v^{n_v}y_v)e^{-\pi \gamma^2 x^2}
$$
Next one has 
$$
\langle 1_{\Z_v}(p_v^{n_v}\bullet)\vert 1_{\Z_v}\rangle =\begin{cases} 1& \text{if}\ n_v\geq 0\\ p_v^{n_v}& \text{if}\ n_v<0
\end{cases}
$$
and thus in all cases 
$$
\langle 1_{R}(\gamma\bullet)\vert 1_{R}\rangle =b^{-1}
$$
We get using \eqref{scalrule}, with $\lambda =e^{-t}$,
$$
\langle \theta(\gamma)\xi\vert \exp(-it\scal)\xi \rangle=\sqrt 2\ \lambda^{1/2}b^{-1}\langle e^{-\pi \gamma^2 x^2}\vert e^{-\pi \lambda^2 x^2}\rangle
$$
$$
\lambda^{1/2}\langle e^{-\pi \gamma^2 x^2}\vert e^{-\pi \lambda^2 x^2}\rangle=\lambda^{1/2}\int_\R e^{-\pi \gamma^2 x^2-\pi \lambda^2 x^2}dx=\frac{\lambda^{1/2}}{\sqrt{\gamma^2+\lambda^2}}=\frac{1}{\sqrt{\lambda^{-1}\gamma^2+\lambda}}
$$
Finally 
$$
b^{-1}\frac{1}{\sqrt{\lambda^{-1}\gamma^2+\lambda}}=\frac{1}{\sqrt{a^2e^{t}+b^2e^{-t}}}
$$
so we obtain the  equality \eqref{expth}. \endproof

\section{The case of a single prime, Lambert series}\label{lambert}
Let $p$ be a prime and $S=\{p,\infty\}$. By \eqref{expth} one gets  
$$
2^{-1/2}\langle\exp(it\scal)\xi_S\vert \xi_S\rangle=\Sigma_p(t)+\Sigma_p(-t)+\frac{1}{\sqrt{e^{-t}+e^{t}}}
$$
where
$$
\Sigma_p(t)=\sum_{n> 0} \frac{1}{\sqrt{p^{2n}e^t+e^{-t}}}=e^{-t/2}\sum_{n> 0}  p^{-n}\frac{1}{\sqrt{1+ p^{-2n}e^{-2t}}}
$$
We can then use the Taylor expansion 
\begin{equation}\label{alphal}
\frac{1}{\sqrt{1+x^2}}=\sum \alpha_\ell x^{2\ell}, \ \   \ \alpha_\ell=(- 4)^{-\ell}\binom{2\ell}{\ell}
\end{equation} 
as well as the summation 
$$
\sum_1^\infty p^{-n-2\ell n}=\frac{p^{-1-2\ell}}{1-p^{-1-2\ell}}
$$
to obtain
\begin{equation}\label{sigmat}
\Sigma_p(t)=e^{-t/2}\sum_0^\infty \alpha_\ell\ e^{-2\ell t}\ \frac{p^{-1-2\ell}}{1-p^{-1-2\ell}}
\end{equation}
One advantage of this formula is that it separates the contribution of the place $\infty$ from the contribution  of the prime $p$ and the latter  tends to $0$ when $p\to \infty$. This suggests that one should look for a $q$-expansion of the result in terms of $q:=p^{-1}$. In fact, one defines in general the Lambert series generating function by the formula 
\begin{equation}\label{Lambert}
L_f(q):=\sum_{n\geq 1} f(n)\frac{q^n}{1-q^n}=\sum_{m\geq 1}(f*1)(m)q^m
\end{equation} 
Here we use the Dirichlet convolution of $f$ and $g, f * g(n):=\sum_{d \mid n} f(d) g(n / d)$, for $n \geq 1$. For $g=1$ this means $f * 1(n):=\sum_{d \mid n} f(d) $.\newline
Taking  into account the prime $p$ in the un-normalized moments \eqref{unmoments}, we define
\begin{equation}\label{normoments}
c(2k,p):=(-1)^k\partial_x^{2k}\left(\frac{1}{\sqrt{e^{-x}+e^{x}}}+\Sigma_p(x)+\Sigma_p(-x)\right)_{x=0}
\end{equation}
\begin{prop}\label{expth3} Let $S=\{p,\infty\}$ where $p$ is a prime and let $c(2k,p)$ as in \eqref{normoments}. One has 
\begin{equation}\label{normoments1}
 c(2k,p)=c_0(2k)+\ (-1)^kL_{f_k}(q), \ q=\frac 1p, \ \ f_k(2\ell+1)=2\, \left(\frac 12+2\ell\right)^{2k}\alpha_l, \  f_k(2\ell)=0.
\end{equation}	
where the Lambert series $L_{f_k}(q)$ is geometrically convergent.
\end{prop}
\proof 
By \eqref{sigmat} one has 
$$
\Sigma_p(t)=\sum_0^\infty \alpha_\ell\ e^{-(\frac 12+2\ell) t}\ \frac{q^{1+2\ell}}{1-q^{1+2\ell}}, \ \ 
$$
$$
\partial_x^{2k}\left(\Sigma_p(x)+\Sigma_p(-x)\right)_{x=0}=2\sum_0^\infty \alpha_\ell\ (\frac 12+2\ell)^{2k}\ \frac{q^{1+2\ell}}{1-q^{1+2\ell}}
$$
which gives the required result using \eqref{Lambert} and \eqref{normoments}. The geometric convergence is ensured by the estimate $\alpha_n\sim (\pi n)^{-1/2}$ while $q=1/p\leq 1/2$.
\endproof
One can now compare numerically the result with the integrals
$$
j(2k,p):=\int_{-\infty}^\infty\vert\frac{1}{1- p^{-\frac{1}{2}-is}}\Gamma(\frac 14-\frac 12 is)\vert^2 s^{2k}ds
$$
One compares the quotients $c(2k,p)/c(0,p)$  with the quotients $j(2k,p)/j(0,p)$. \newline
For $p=2$ the numerical computation gives the following values for $j(2k,p)/j(0,p)$
$$
\{0.157505,0.202341,1.479,43.9807,2503.76,178096.,1.40009\times 10^7,1.19254\times 10^9\}
$$
One recovers the following values for $c(2k,p)/c(0,p)$ using the first $40$ terms of the Lambert series 
$$
\{0.157505,0.202341,1.479,43.9807,2503.76,178096.,1.40009\times 10^7,1.19257\times 10^9\}
$$
\subsection{Observed integrality of the determinant $D_n(q)$}\label{sectq}

The next goal is to determine the deformation of the coefficients $a(n)$ due to the presence of a single prime. 
In fact one observes the following key fact
\begin{fact}\label{keyfact} The $q$-expansion of the determinant $D_n(q)$ is of the form, 
\begin{equation}\label{qexpandn}
D_n(q)=D_n(0)\left(\sum_0^\infty \mu_{n,k}\ q^k\right), \  \ \mu_{n,k}\in \Z[\sqrt 2,\frac 12]=\Z[\frac{1}{\sqrt 2}]
\end{equation}	
\end{fact}
This fact is quite meaningful because the terms $D_n(0)$ are of the form for $0\leq n\leq 10$ 
$$
\frac{1}{\sqrt{2}},\frac{1}{4},\frac{3}{8 \sqrt{2}},\frac{135}{64},\frac{42525}{128 \sqrt{2}},\frac{602791875}{1024},\frac{281970969328125}{4096 \sqrt{2}},\frac{12002806286149194140625}{65536},$$ $$\frac{7663946349837194629663623046875}{131072 \sqrt{2}},\frac{748709847905086347512925275890191650390625}{1048576},$$ $$\frac{6948614854850430587049622200362718103307453155517578125}{4194304 \sqrt{2}}
$$
where the denominators are powers of $\sqrt 2$ and the numerators are prime to $2$. 
We let 
\begin{equation}\label{qexpandn1}
\Delta_n(q)=\sum_0^\infty \mu_{n,k}\ q^k =D_n(q)/D_n(0)
\end{equation}
The first instances of these $q$-series are as follows
$$
\Delta_1(q)=1+\sqrt{2}\, q+(\sqrt{2}-4)q^2+(\frac{25}{\sqrt{2}}-8)q^3+(\sqrt{2}+40)q^4+(36-\frac{229 }{4 \sqrt{2}})q^5+\ldots$$
$$
\Delta_2(q)=1+\frac{7}{2 \sqrt{2}}q+\frac{1}{4} \left(7 \sqrt{2}-22\right)q^2-(11+\frac{137}{4 \sqrt{2}})q^3+(19+\frac{7}{2 \sqrt{2}})q^4+$$ $$+\frac{3}{32} \left(3023 \sqrt{2}+144\right)q^5+\frac{1}{8} \left(1399 \sqrt{2}-4604\right)q^6+\ldots
$$
$$
\Delta_3(q)=1+\frac{9}{4 \sqrt{2}}q+\frac{1}{16} \left(18 \sqrt{2}-143\right)q^2+\frac{11}{16} \left(51 \sqrt{2}-26\right)q^3+\frac{1}{8} \left(9 \sqrt{2}+2387\right)q^4$$ $$+\frac{1}{64} \left(18524-33837 \sqrt{2}\right)q^5+\frac{1}{16} \left(2033 \sqrt{2}-111923\right)q^6+\ldots
$$
$$
\Delta_4(q)=1+\frac{107}{32 \sqrt{2}}q+\frac{1}{256} \left(428 \sqrt{2}-2743\right)q^2+\frac{1}{128} \left(-5557 \sqrt{2}-2743\right)q^3+$$ $$+\frac{1}{128} \left(214 \sqrt{2}+40327\right)q^4+\frac{17}{512} \left(56537 \sqrt{2}+9166\right)q^5+\ldots
$$
$$
\Delta_5(q)=1+\frac{151}{64 \sqrt{2}}q+\frac{4832 \sqrt{2}-57181}{4096}q^2+\frac{132344 \sqrt{2}-57181}{2048}q^3+$$ $$+\frac{2416 \sqrt{2}+2015041}{2048}q^4+\frac{3972901-9954328 \sqrt{2}}{4096}q^5+\ldots 
$$
$$
\Delta_6(q)=1+\frac{835}{256 \sqrt{2}}q+\frac{53440 \sqrt{2}-520049}{32768}q^2-\frac{1243856 \sqrt{2}+520049}{16384}q^3+$$ $$+\frac{26720 \sqrt{2}+18184313}{16384}q^4+\frac{405931619 \sqrt{2}+71697154}{65536}q^5+\ldots $$
$$
\Delta_7(q)=1+\frac{1241}{512 \sqrt{2}}q+\frac{317696 \sqrt{2}-4983643}{262144}q^2+\frac{7 \left(1869248 \sqrt{2}-711949\right)}{131072}q^3+$$ $$+\frac{158848 \sqrt{2}+302791807}{131072}q^4+\frac{2402399884-7326921115 \sqrt{2}}{1048576}q^5+\ldots 
$$
$$
\Delta_8(q)=1+\frac{26291}{8192 \sqrt{2}}q+\frac{13460992 \sqrt{2}-176188271}{8388608}q^2-\frac{476628608 \sqrt{2}+176188271}{4194304}q^3+$$ $$\frac{6730496 \sqrt{2}+11003551439}{4194304}q^4+\frac{174647316856 \sqrt{2}+2004033963205}{67108864 \sqrt{2}}q^5+\ldots
$$
All these $q$-series have coefficients in the ring $\cR=\Z[\frac{1}{\sqrt 2}]$ and leading coefficient $1$. Such $q$-series form a commutative group  $\widehat{W}(\cR)$ under multiplication.
\subsection{Deformation of the Jacobi coefficients}
\begin{cor}\label{ansq1} Let $a(n)^2(q)$ as in \eqref{ansq}. One has $a(n)^2(q)/a(n)^2(0)\in \widehat{W}(\cR)$
$$
a(n)^2(q)=(n+\frac 12)(n+1)\frac{\Delta_{n-1}(q)\Delta_{n+1}(q)}{\Delta_n(q)^2},\ \ 
$$	
\end{cor}
\proof This follows directly from \eqref{ansq} and \eqref{qexpandn1} using $a(n)^2(0)=(n+\frac 12)(n+1)$.\endproof 
In fact the additive group $\widehat{W}(\cR)$ underlies a ring structure which comes from the Almkvist Witt ring and which is in fact a $\lambda$-ring of Witt vectors.  With these structures one can think of the ratio $\frac{\Delta_{n-1}(q)\Delta_{n+1}(q)}{\Delta_n(q)^2}$ as a second finite difference, since product is interpreted as a sum. In order to look for a relevance of this structure one takes the ghost components, which gives the following matrix whose lines give the first ghost components ${\rm gh}_j(D_k)$, $1\leq j\leq 4$.
$$
\left(
\begin{array}{cccc}
 -\sqrt{2} & -2 \left(\sqrt{2}-5\right) & \frac{1}{2} \left(60-103 \sqrt{2}\right) & -12 \left(5 \sqrt{2}-1\right) \\
 -\frac{7}{2 \sqrt{2}} & \frac{1}{8} \left(137-28 \sqrt{2}\right) & \frac{1}{32} \left(377 \sqrt{2}+1644\right) & \frac{1}{64} (-9) \left(1176 \sqrt{2}+503\right) \\
 -\frac{9}{4 \sqrt{2}} & \frac{1}{32} \left(653-72 \sqrt{2}\right) & \frac{1}{256} (-3) \left(11793 \sqrt{2}-5224\right) & \frac{-139824 \sqrt{2}-631007}{1024} \\
 -\frac{107}{32 \sqrt{2}} & \frac{55337-6848 \sqrt{2}}{2048} & \frac{8802037 \sqrt{2}+10624704}{131072} & \frac{-1086489984 \sqrt{2}-5575217647}{4194304} \\
 -\frac{151}{64 \sqrt{2}} & \frac{251525-19328 \sqrt{2}}{8192} & \frac{96585600-258529321 \sqrt{2}}{1048576} & \frac{-14460397824 \sqrt{2}-185707273519}{67108864} \\
 -\frac{835}{256 \sqrt{2}} & \frac{4857617-427520 \sqrt{2}}{131072} & \frac{9491428673 \sqrt{2}+7461299712}{67108864} & -\frac{5 \left(1208835886080 \sqrt{2}+15632499485747\right)}{17179869184} \\
\end{array}
\right)
$$
The ghost components of a series $f(q)$ are given by the coefficients of the $q$-expansion of the  $-\frac{q f'(q)}{f(q)}$. They define characters of the ring  $\widehat{W}(\cR)$.\newline

\section{Proof of integrality} \label{secinteg}
We now try to understand the meaning of the divisibility property expressed by Fact \ref{keyfact}. Except for powers of $2$ it is divisibility of all terms of the $q$-expansion by the product
\begin{equation}\label{psifun}
\Psi(n):=\prod _{k=1}^n (2 k)!
\end{equation}
In fact with the normalization chosen for the moments one gets the equality 
\begin{equation}\label{psifun1}
D_n(0)=2^{-n^2-\frac{3 n}{2}-\frac{1}{2}}\, \Psi(n)
\end{equation}
The first attempt to understand Fact \ref{keyfact} is to view the determinant $D_n$ as a perturbation of $D_n(0)$ where the matrix entries are perturbed using $c(2j)\mapsto c(2j)+\epsilon_j$. One can even take all $\epsilon_j=0$ except $\epsilon_1=1$. Then one obtains the following values for the quotient of the perturbed determinant $D'_n$ divided by $D_n(0)$, 
$$
-\frac{4 \sqrt{2}}{3}-3,\frac{1}{135} \left(-229 \sqrt{2}-601\right),\frac{-599287 \sqrt{2}-786735}{85050}$$ $$\frac{-914620363 \sqrt{2}-1191630600}{107163000},\frac{-134998296477451 \sqrt{2}-143229230027100}{8911675080000}
$$
Thus it is not true that the divisibility Fact \ref{keyfact} comes simply from perturbation. The next idea is that it could come from a matrix factorization. Thus one computes the inverse of the Hankel matrices $\cD_n(q)$ evaluated at $q=0$. One gets matrices $\cD_n(0)^{-1}$ of the form
$$
\cD_5(0)^{-1}=\frac{1}{\sqrt 2}\left(
\begin{array}{cccccc}
 \frac{211}{84} & 0 & -\frac{25}{21} & 0 & \frac{1}{21} & 0 \\
 0 & \frac{219661}{28350} & 0 & -\frac{3874}{2835} & 0 & \frac{422}{14175} \\
 -\frac{25}{21} & 0 & \frac{904}{315} & 0 & -\frac{44}{315} & 0 \\
 0 & -\frac{3874}{2835} & 0 & \frac{1504}{2835} & 0 & -\frac{8}{567} \\
 \frac{1}{21} & 0 & -\frac{44}{315} & 0 & \frac{4}{315} & 0 \\
 0 & \frac{422}{14175} & 0 & -\frac{8}{567} & 0 & \frac{8}{14175} \\
\end{array}
\right)
$$
One then multiplies this matrix by the Hankel matrix depending on $q$, and one obtains a matrix of functions of $q$. But even for $n=3$ the obtained matrix does not have coefficients in $\Z[\frac{1}{\sqrt 2}]$ since it contains the diagonal term $-\frac{19 q^4}{6 \sqrt{2}}$ but its trace seems to be in $\Z[\frac{1}{\sqrt 2}]$ which suggests that a suitable similar matrix could well have all its coefficients in $\Z[\frac{1}{\sqrt 2}]$. \newline
In fact one can start by considering the terms $\mu_{n,1}\, q$ in $\Delta_n(q)$ which are linear in $q$. One obtains the following formula for these terms using the guess given by the computer in terms of the hypergeometric function $\, _2F_1\left(1,n+\frac{3}{2};n+2;-1\right)$. One has in general
$$
{ }_2 F_1(a, b ; c ; z)=\sum_{k=0}^{\infty} \frac{(a)_k(b)_k}{(c)_k} \frac{z^k}{k !}, $$  
$$(a)_k= \begin{cases}1 & \text { if } k=0 \\ a(a+1) \cdots(a+k-1) & \text { if } k>0\end{cases}
$$
Thus one gets, since the term $(a)_k$ for $a=1$ cancells $k!$,
$$
(-1)^{n+1}\left(\frac{\prod _{j=0}^{n} \left(j+\frac{1}{2}\right)}{(n+1)!}\right){}_2F_1\left(1,n+\frac{3}{2};n+2;-1\right)=\sum _{k=n+1}^{\infty } \frac{(-1)^k \left(\prod _{j=0}^{k-1} \left(j+\frac{1}{2}\right)\right)}{k!}
$$
while one has
$$
\sum _{k=1}^{\infty } \frac{(-1)^k \left(\prod _{j=0}^{k-1} \left(j+\frac{1}{2}\right)\right)}{k!}=\frac{1}{2} \left(\sqrt{2}-2\right)
$$
\subsection{The  $q$-matrices $\cA_n(q)$}
This suggests that with $\alpha_k$ as above in \eqref{alphal}, one should have
\begin{equation}\label{mun1}
\mu_{n,1}=2 \sqrt{2}\ \sum_0^n\alpha_k\in \Z[\frac{1}{\sqrt 2}]
\end{equation}
One needs to prove this fact directly and we use
\begin{lem}\label{detastrace} $(i)$~Let $\cE_n(q):=\cD_n(q)-\cD_n(0)$ and $\cA_n(q):=\cD_n(0)^{-1}\cE_n(q)$. One has 
\begin{equation}\label{Ln0}
\det(\cD_n(q))=\det(\cD_n(0))\sum \Tr(\wedge^m \cA_n(q))	
\end{equation}
$(ii)$~The coefficient	$\mu_{n,1}$ is the trace of the matrix which is the product 
\begin{equation}\label{Ln}
L_n:=\cD_n(0)^{-1}(\partial_q \cD_n)(0)	
\end{equation}
\end{lem}
\proof $(i)$~Follows from the equalities 
$$
\cD_n(q)=\cD_n(0)\left(1+\cA_n(q)\right), \ \ \det(1+\cA)=\sum \Tr(\wedge^m \cA)
$$  
$(ii)$~Follows from $(i)$. \endproof  
The matrix $\partial_q \cD_n(0)$ has a simple form as a Hankel matrix with entries signed powers of $2$
$$
\partial_q \cD_5(0)=\left(
\begin{array}{cccccc}
 2 & 0 & -\frac{1}{2} & 0 & \frac{1}{8} & 0 \\
 0 & -\frac{1}{2} & 0 & \frac{1}{8} & 0 & -\frac{1}{32} \\
 -\frac{1}{2} & 0 & \frac{1}{8} & 0 & -\frac{1}{32} & 0 \\
 0 & \frac{1}{8} & 0 & -\frac{1}{32} & 0 & \frac{1}{128} \\
 \frac{1}{8} & 0 & -\frac{1}{32} & 0 & \frac{1}{128} & 0 \\
 0 & -\frac{1}{32} & 0 & \frac{1}{128} & 0 & -\frac{1}{512} \\
\end{array}
\right)$$
The matrix $L_n$ of \eqref{Ln} is simpler than $\cD_n(0)^{-1}$ and takes the following form
$$
L_5=\left(
\begin{array}{cccccc}
 \frac{45}{8 \sqrt{2}} & 0 & -\frac{45}{32 \sqrt{2}} & 0 & \frac{45}{128 \sqrt{2}} & 0 \\
 0 & -\frac{971}{240 \sqrt{2}} & 0 & \frac{971}{960 \sqrt{2}} & 0 & -\frac{971}{3840 \sqrt{2}} \\
 -\frac{23}{6 \sqrt{2}} & 0 & \frac{23}{24 \sqrt{2}} & 0 & -\frac{23}{96 \sqrt{2}} & 0 \\
 0 & \frac{3}{4 \sqrt{2}} & 0 & -\frac{3}{16 \sqrt{2}} & 0 & \frac{3}{64 \sqrt{2}} \\
 \frac{1}{6 \sqrt{2}} & 0 & -\frac{1}{24 \sqrt{2}} & 0 & \frac{1}{96 \sqrt{2}} & 0 \\
 0 & -\frac{1}{60 \sqrt{2}} & 0 & \frac{1}{240 \sqrt{2}} & 0 & -\frac{1}{960 \sqrt{2}} \\
\end{array}
\right)
$$
where one can see that the lines are determined by their first non-zero term by multiplication by powers of $-\frac 14$. This is proved since the columns of $\partial_q \cD_n(0)$ are obtained from the first two columns by multiplication by powers of $-\frac 14$. In fact one needs to analyze the matrix $\partial_q \cD_n(0)$ which is of rank $2$ and to write it as a sum of two multiples of rank one projections. The eigenvalues are 
$$
\lambda_1(n)=\frac{2}{15} 16^{1-m} \left(16^m-1\right),\  m=\left\lceil \frac{n+1}{2}\right\rceil, \ \ \lambda_2(n)=-\frac{2}{15} 4^{1-2 m} \left(4^{2 m}-1\right), \ m=\left\lceil \frac{n}{2}\right\rceil
$$
The eigenvectors are of the following form (for the two eigenvalues)
$$
\xi_1(n)=\left(1,0,-\frac{1}{4},0,\frac{1}{16},0,\ldots\right), \ \ \xi_2(n)=\left(0, -\frac i2, 0, \frac i8, 0, -\frac{i}{32}\ldots \right)
$$
The next step is to get an explicit formula for the inverse Hankel matrix $\cD_n(0)^{-1}$, one has
\begin{lem} \label{inverseh} One has $\cD_n(0)^{-1}=\sqrt 2\, \cB\, \cB^*$ where $\cB=(b_{k,\ell})_{k\leq \ell}$ and the coefficients $b_{k,\ell}$ are given by $\cP_\ell(x)=\sum b_{k,\ell}\, x^k$ where the polynomials $\cP_m$ are defined by the formula
\begin{equation}\label{basicpols}
\cP_m(x):=\sqrt{(2 m)!}\sum_0^m (-1)^k 2^{3 k-m}\frac{ \prod _{j=0}^{k-1} \left(j-\frac{i x}{2}+\frac{1}{4}\right)}{(2 k)! (m-k)!}
\end{equation}	
\end{lem}
\proof One uses the result of A.C.~Aitken (see \cite{berg}) which asserts that the inverse of the Hankel matrix defined by the normalized moments is equal to $\cB\, \cB^*$ where the matrix entries $b_{k,\ell}$ of the upper triangular matrix $\cB_n$ are given by $\cP_\ell(x)=\sum b_{k,\ell}\, x^k$ where the polynomials $\cP_\ell$ are the orthonormal polynomials for the given normalized moments. The orthonormal polynomials in our case are given by Proposition 3.3 of \cite{c2m}, \ie  \eqref{basicpols}.  Our moments are not normalized and the normalized ones are obtained by multiplication by $\sqrt 2$, this gives the overall factor $\sqrt 2$.\endproof 
Note that, unlike for $\cD_n$, the entries of the  triangular matrices $\cB_n$ do not change with $n$
$$
\cB_5=\left(
\begin{array}{cccccc}
 1 & 0 & \frac{1}{\sqrt{6}} & 0 & \frac{\sqrt{\frac{5}{14}}}{2} & 0 \\
 0 & i \sqrt{2} & 0 & \frac{7 i}{3 \sqrt{5}} & 0 & \frac{211 i}{90 \sqrt{7}} \\
 0 & 0 & -\sqrt{\frac{2}{3}} & 0 & -\frac{1}{3} \left(11 \sqrt{\frac{2}{35}}\right) & 0 \\
 0 & 0 & 0 & -\frac{2 i}{3 \sqrt{5}} & 0 & -\frac{10 i}{9 \sqrt{7}} \\
 0 & 0 & 0 & 0 & \frac{\sqrt{\frac{2}{35}}}{3} & 0 \\
 0 & 0 & 0 & 0 & 0 & \frac{2 i}{45 \sqrt{7}} \\
\end{array}
\right)
$$
We then obtain the following formula for the coordinates of the vectors $\cB^*(\xi_j(n))$,
\begin{equation}\label{vectorsj}
\left(\cB^*(\xi_1(n))\right)_\ell=\frac 12\left(\cP_\ell(\frac i2)+\cP_\ell(-\frac i2)\right),\ \ \left(\cB^*(\xi_2(n))\right)_\ell=\frac 12 \left(\cP_\ell(\frac i2)-\cP_\ell(-\frac i2)\right)
\end{equation}
Now one has, using \eqref{basicpols}, the equalities 
\begin{equation}\label{basicpols1}
\cP_m(i/2)=\sqrt{(2 m)!}\sum_0^m (-1)^k 2^{3 k-m}\frac{ \prod _{j=0}^{k-1} \left(j+\frac{1}{2}\right)}{(2 k)! (m-k)!}
\end{equation}
\begin{equation}\label{basicpols2}
\cP_m(-i/2)=\sqrt{(2 m)!}\sum_0^m (-1)^k 2^{3 k-m}\frac{ \prod _{j=0}^{k-1} j}{(2 k)! (m-k)!}
\end{equation}
\begin{lem} \label{tracelem} The trace $\Tr(L_n)$ of $L_n$ is given by
\begin{equation}\label{traceform}
\frac {\sqrt 2}{4}\lambda_1(n)\Vert \xi_1(n)\Vert^{-2}\sum \left(\cP_\ell(\frac i2)+\cP_\ell(-\frac i2)\right)^2+\frac {\sqrt 2}{4}\lambda_2(n)\Vert \xi_2(n)\Vert^{-2}\sum \left(\cP_\ell(\frac i2)-\cP_\ell(-\frac i2)\right)^2
\end{equation}	
\end{lem}
 \proof This follows from the equality for the rank two matrix $\partial_q \cD_n(0)$,
 $$
 \partial_q \cD_n(0)=\lambda_1(n)\Vert \xi_1(n)\Vert^{-2}\vert \xi_1(n)\rangle \langle \xi_1(n)\vert+\lambda_2(n)\Vert \xi_2(n)\Vert^{-2}\vert \xi_2(n)\rangle \langle \xi_2(n)\vert
 $$
and from Lemma \ref{inverseh} together with \eqref{vectorsj}. \endproof 
Next, note that $\lambda_1(n)\Vert \xi_1(n)\Vert^{-2}=2$ and that $\lambda_2(n)\Vert \xi_2(n)\Vert^{-2}=-2$, so that we obtain  
\begin{equation}\label{traceform1}
\Tr(L_n)=2\sqrt 2\sum \cP_\ell(i/2)\cP_\ell(-i/2)
\end{equation}	
One now has the following computation of $\cP_\ell(i/2)\cP_\ell(-i/2)$
\begin{prop}	
 \label{tracelem1} $(i)$~One has the equalities for any $m \in \N$,
\begin{equation}\label{traceform2}
 \cP_m(i/2)=(-1)^m 2^{-m} \frac{\sqrt{(2 m)!}}{m!},\  \cP_m(-i/2)=2^{-m} \frac{\sqrt{(2 m)!}}{m!}\end{equation}
\begin{equation}\label{traceform3}
  \cP_\ell(i/2)\cP_\ell(-i/2)=\alpha_\ell
\end{equation}	
$(ii)$~Equality \eqref{mun1} holds \ie $\mu_{n,1}=2 \sqrt{2}\ \sum_0^n\alpha_k\in \Z[\frac{1}{\sqrt 2}]$ for all $n\in \N$.
\end{prop}
\proof $(i)$~One has $\cP_m(i/2)=\cP_m(-i/2)$ for $m$ even since in that case the polynomial $\cP_m$ is real and by \eqref{basicpols1} the value $\cP_m(i/2)$ is real. For $m$ odd the polynomial $\cP_m$ is purely imaginary but the value $\cP_m(i/2)$ is real still, so one gets $\cP_m(i/2)=-\cP_m(-i/2)$ for $m$ odd. One has by \eqref{basicpols2} that all the terms in the sum involved for $\cP_m(-i/2)$ are equal to $0$ except the first term. Equivalently one can use \eqref{basicpols} and thus one gets \eqref{traceform2}.
 Then \eqref{traceform3} follows.\newline 
 $(ii)$~Follows from \eqref{traceform1}. \endproof

 One can moreover give a formula for the entry $L_n(1,1)$ using the Pochhammer symbol $(a)_{n}:=a(a+1)\ldots (a+n-1)$. One gets the formula 
 $$
 L_n(1,1)=\frac{5 \left(\frac{9}{4}\right)_{n-1}}{(2)_{n-1}}
 $$
 \subsection{Higher order terms in the $q$-expansion of the determinant $D_n(q)$}
 We now  use Lemma \ref{detastrace} in order to understand the higher order terms in the $q$-expansion of the determinant $D_n(q)$. We describe a matrix similar to $\cA_n(q):=\cD_n(0)^{-1}\cE_n(q)$ and which hence gives the same result for the trace on the wedge powers. 
 \begin{prop} \label{similar} The matrix $\cA_n(q):=\cD_n(0)^{-1}\cE_n(q)$  is similar (by a conjugation by an invertible matrix) to the direct sum of two matrices $\cA_n^+(q)\oplus \cA_n^-(q)$.  Each $\cA_n^\pm(q)$ is a Lambert series  
 	\begin{equation}\label{lambertrank1}
\cA_n^\pm(q)=\sum_{\ell\geq 0} \alpha_\ell \frac{q^{2\ell+1}}{1-q^{2\ell+1}}T_\ell^\pm
\end{equation}	
where  $T_\ell^\pm$ are the rank one operators $T_\ell^\pm:=\pm 2^{3/2}\vert \eta_\ell^\pm\rangle \langle \eta_\ell^\pm\vert$ with the vectors  $\eta_\ell^\pm$ with components obtained from the  sequence of real numbers $\eta_\ell(m):=\cP_m(-i(\frac 12+2\ell))$ by
\begin{equation}\label{tlplusminus}
 \eta_\ell^\pm(k):= \frac{1\pm(-1)^k}{2}\eta_\ell(k) 
\end{equation}
 \end{prop}
 \proof By Lemma \ref{detastrace} $(i)$, one has $\cA_n(q):=\cD_n(0)^{-1}\cE_n(q)$ and, by  Lemma \ref{inverseh}, one has the equality $\cD_n(0)^{-1}=\sqrt 2\cB\, \cB^*$ where $\cB=(b_{k,\ell})_{k\leq \ell}$ while the coefficients $b_{k,\ell}$ are given by $\cP_\ell(x)=\sum b_{k,\ell}\, x^k$ where the $\cP_\ell$ are given in \eqref{basicpols}. Moreover  the polynomials $\cP_m$ are even for $m$ even and odd for $m$ odd. The triangular matrix $\cB$ is invertible since each diagonal term is the term of  degree $\ell$ in $\cP_\ell$ which is non-zero since $\cP_\ell$  has degree $\ell$. Thus the matrix $\cA_n(q)$ is similar to $\cA'_n(q)$ where
 \begin{equation}\label{bstarb}
 \cA'_n(q):=\cB^{-1}\cA_n(q)\cB=\sqrt 2 \ \cB^*\cE_n(q)\cB
 \end{equation}
 By Proposition \ref{expth3}, $\cE_n(q)$ is the Lambert series 
\begin{equation}\label{lambplusminus}
 \cE_n(q)=\sum_{\ell\geq 0} \alpha_\ell \frac{q^{2\ell+1}}{1-q^{2\ell+1}}E_{n,\ell}
 \end{equation} 
 where for each $\ell\geq 0$, the matrix $E_{n,\ell}$ is the   Hankel matrix associated with the even moments given by  \eqref{normoments1}, \ie $$\epsilon(2k)=2(-1)^k\left(\frac 12 +2\ell\right)^{2k}$$
 Now for fixed $n$ the vector space $\cH_n$ of dimension $n+1$ of vectors with coordinates $\xi_j$, $j\in \{0,\ldots,n\}$, decomposes as the direct sum of the two subspaces $\cH_n^\pm$ of vectors whose all odd coordinates (resp. even) vanish. We use the inner product on $\cH$ given by the norm square $\sum \vert \xi_j\vert^2$. Given $\xi,\eta\in \cH$ the rank one operator $T=\vert \xi\rangle\langle \eta\vert$ is defined by the equality  $T(\alpha):=\langle \eta\vert \alpha\rangle \xi$ where the sesquilinear form $\langle \eta\vert \alpha\rangle$ is linear in $\alpha$ and antilinear in $\eta$. The entries of the matrix $T=\vert \xi\rangle\langle \eta\vert$ are the $T_{i,j}=\xi_i\overline{\eta_j}$. The Hankel matrix $H(c)$ associated to any even moments $c(2k)$ preserves the subspaces $\cH_n^\pm$, and hence its determinant factorizes as the product of the determinants of the restrictions $H^\pm(c)$
 $$
 \det(H(c))=\det(H^+(c))\det(H^-(c))
 $$
 The Hankel matrix $E_{n,\ell}$ associated to the even moments $\epsilon(2k)$ is the direct sum $E_{n,\ell}=E^+_{n,\ell}\oplus E^-
 _{n,\ell}$ and the matrices $E^\pm_{n,\ell}=\pm2\vert \gamma^\pm_\ell\rangle\langle \gamma^\pm_\ell\vert$ are of rank one, where the vectors  $\gamma^\pm_\ell\in \cH_n^\pm$ are
 $$
\gamma^+_\ell=\left(1,0,-(\frac{1}{2}+2\ell)^2,0,(\frac{1}{2}+2\ell)^4,0,\ldots\right), \ \ \gamma^-_\ell=\left(0, (\frac{1}{2}+2\ell)i, 0, (\frac{1}{2}+2\ell)^3i, 0, \ldots \right)
$$
 which correspond to the decomposition in even and odd parts of the vector $\gamma_\ell\in \cH_n$ whose components are the powers of $(\frac{1}{2}+2\ell)i$. Thus one has 
 \begin{equation}\label{gamplusminus}
 \gamma_\ell^\pm(k):= \frac{1\pm(-1)^k}{2}z^k, \ \ z=(\frac{1}{2}+2\ell)i. 
\end{equation}
 The matrix of $\cB^*$ is the lower triangular matrix with entries $(\cB^*)_{m,k}=\overline{b_{k,m}}$ where  $\cP_m(x)=\sum b_{k,m}\, x^k$ with the polynomials  defined by the formula \eqref{basicpols}. The coefficients $b_{k,m}$ are zero when $m-k$ is odd, they are real for $m$ even and purely imaginary for $m$ odd. Thus one gets
 $$
 \sum (\cB^*)_{m,k}\gamma_\ell^+(k)=\cP_m(z)\ \ \text{if}\ \ m \ \text{is even}, \ 0 \ \text{otherwise}
 $$
Similarly, using $\overline{b_{k,m}}=-b_{k,m}$ for $m$ and $k$ odd, one gets
 $$
 \sum (\cB^*)_{m,k}\gamma_\ell^-(k)=-\cP_m(z)=\cP_m(-z)\ \ \text{if}\ \ m \ \text{is odd}, \ 0 \ \text{otherwise}
 $$
 One has $\cB^*\vert \xi\rangle\langle \eta\vert\cB=\vert\cB^*( \xi)\rangle\langle \cB(\eta)\vert$ for any vectors $\xi, \eta$ and thus, using \eqref{bstarb}, \eqref{lambplusminus} we obtain
 $$
 \cA'_n(q)=\sqrt 2 \ \cB^*\cE_n(q)\cB=\sqrt 2 \ \sum_{\ell\geq 0} \alpha_\ell \frac{q^{2\ell+1}}{1-q^{2\ell+1}}\cB^*E_{n,\ell}\cB
 $$
 $$
 \cB^*E_{n,\ell}\cB=\cB^*E^+_{n,\ell}\cB\oplus \cB^*E^-_{n,\ell}\cB, \ \  
 \cB^*E^\pm_{n,\ell}\cB=\pm 2\vert \eta_\ell^\pm\rangle \langle \eta_\ell^\pm\vert
 $$
 where the vectors $\eta_\ell^\pm\in \cH_n^\pm$ are the even and odd parts of the vector $\eta_\ell$ with components $\cP_m(-z)$ which are real numbers explicitly given using the formula \eqref{basicpols} as
 \begin{equation}\label{pmz}
\cP_m(-z)=\sqrt{(2 m)!}\sum_0^m (-1)^k 2^{3 k-m}\frac{ \prod _{j=0}^{k-1} \left(j-\ell\right)}{(2 k)! (m-k)!}
\end{equation}
 We thus obtain  
 \begin{equation}\label{lambertrank1g}
\cA_n^\pm(q)=\pm 2\sqrt 2 \sum_{\ell\geq 0} \alpha_\ell \frac{q^{2\ell+1}}{1-q^{2\ell+1}}\vert \eta_\ell^\pm\rangle \langle \eta_\ell^\pm\vert
\end{equation}
 which is the required decomposition \eqref{lambertrank1}.\endproof 
 Note that \eqref{pmz}, which gives the components of the vector $\eta_\ell$, can be written as  
 \begin{equation}\label{pmzg}
\eta_\ell(m)=2^{-m}\sqrt{(2 m)!}\sum_{k=0}^\ell  2^{3 k}\frac{ \prod _{j=0}^{k-1} \left(\ell-j\right)}{(2 k)! (m-k)!}
\end{equation}
 since  $\prod _{j=0}^{k-1} \left(\ell-j\right)=0$ for $k>\ell$. In turns this can be rewritten as 
 \begin{equation}\label{pmzg1}
\eta_\ell(m)=\frac{2^{-m}\sqrt{(2 m)!}}{m!}\sum_{k=0}^\ell  2^{3 k}\frac{k! \prod _{j=0}^{k-1} \left(\ell-j\right)}{(2 k)! }\binom{m}{k}
\end{equation}
Thus the dependence on $m$ splits as a product of the term $\frac{2^{-m}\sqrt{(2 m)!}}{m!}$ whose square is $(-1)^m\alpha_m$ (as in \eqref{alphal}) and of the term 
 \begin{equation}\label{pmzg2}
p_\ell(m):=\sum_{k=0}^\ell  2^{3 k}\frac{k! \prod _{j=0}^{k-1} \left(\ell-j\right)}{(2 k)! }\binom{m}{k}=\sum_{k=0}^\ell  2^{3 k}\binom{m}{k}\binom{\ell}{k}\binom{2k}{k}^{-1}
\end{equation}
which is a polynomial of degree $\ell$ as a function of $m$ and is a symmetric function of $m$ and $\ell$.
 \begin{cor}\label{factorization} The determinant $D_n(q)$ factorizes as the product 
 $D^+_n(q)\times D^-_n(q)$ where 
 \begin{equation}\label{factori1}
 D^\pm_n(q)=D^\pm_n(0)\sum \Tr(\wedge^m \cA^\pm_n(q))	
 \end{equation}
Moreover one has $D^+_n(0)=\prod _{k=0}^n 2^{-4 k-\frac{1}{2}} (2 (2 k))!$ and $D_n^-(0)=\prod _{k=0}^n 2^{-4 k-\frac{5}{2}} (2 (2 k+1))!$.
 \end{cor}
\proof The factorization \eqref{factori1} follows the decomposition of matrices into even and odd parts, which gives a factorization of their determinants. One then  uses Lemma \ref{detastrace} applied to the decomposition $\cA_n^+(q)\oplus \cA_n^-(q)$. To compute $D^\pm_n(0)$ one notices that when $n=2k$ is even the dimensions of $\cH_n^\pm$ are respectively $k+1$ and $k$ while for $n=2k+1$ the dimensions of $\cH_n^\pm$ are both equal to $k+1$. Moreover when one passes from $n$ to $n+1$ the matrix does not change in one of the two spaces, for $n$ even the matrix  $\cA_n^+(q)$ does not change and for $n$ odd the matrix $\cA_n^+(q)$. Since one knows the ratio $D_{n+1}(0)/D_{n}$ one obtains the required formulas. \endproof 
\subsection{Integrality and Catalan numbers}
In fact the computations show that there is a refined form of Fact \ref{keyfact}, namely that it seems to hold for each $D^\pm_n(q)/D^\pm_n(0)$. In order to prove these facts one investigates the entries of the matrices $\cA_n^\pm(q)$ using Proposition \ref{similar}, and one finds that they belong to the ring $\ccyc[\frac 12]$ where $\ccyc$ is the ring of cyclotomic integers. In order to prove this, one needs the following integrality result:
\begin{lem}\label{integralitylem}$(i)$~For each $k\leq a$ the central binomial coefficient $\binom{2k}{k}$ divides the product $\binom{2a}{a}\binom{a}{k}$.\newline
$(ii)$~Let $p_\ell(m)$ be given by \eqref{pmzg2}, then for any $\ell,m$ the following number is an integer 
$$
\binom{2\ell}{\ell}\, p_\ell(m)\in \Z
$$
$(iii)$~For any $\ell,m$ the following number is an integer
\begin{equation}\label{sigmalm}
	\sigma(\ell,m):=\binom{2\ell}{\ell}\binom{2m}{m}p_\ell(m)^2
\end{equation}
\end{lem}
\proof $(i)$~Let $k=a-u$. One has 
$$
\binom{2a}{a}/\binom{2k}{k}=\frac{(2a)!}{(2(a-u))!}((a-u)!)^2(a!)^{-2}
$$
$$
\frac{(2a)!}{(2(a-u))!}=\prod_{0\leq j\leq 2u-1} (2a-j), \ \ \frac{a!}{(a-u)!}=\prod_{0\leq j\leq u-1} (a-j)
$$
Thus one gets 
$$
\frac{(2a)!}{(2(a-u))!}/\frac{a!}{(a-u)!}=2^u\prod_{0\leq j\leq u-1} (2a-2j-1)
$$
so that 
$$
\binom{2a}{a}/\binom{2k}{k}=2^u\prod_{0\leq j\leq u-1} (2a-2j-1)/\prod_{0\leq j\leq u-1} (a-j)
$$
Moreover 
$$
\binom{a}{k}=\binom{a}{u}=\left(\prod_{0\leq j\leq u-1} (a-j)\right)/u!
$$
Thus one obtains the  following equality 
 \begin{equation}\label{halfinteg}
\binom{a}{k}\binom{2a}{a}/\binom{2k}{k}=\frac{2^u}{u!}\prod_{0\leq j\leq u-1} (2a-2j-1)=2^{2u}\binom{a-\frac 12}{u}=2^{2(a-k)}\binom{a-\frac 12}{a-k}
 \end{equation}
 To obtain that this number is an integer one reduces to the computation of 
 $$
 2^{2k-1}\binom{\frac 12}{k}=(-4)^{k-1}(k!)^{-1}\prod_{0\leq j\leq k-2} \left(\frac 12 +j\right)
 $$
 These numbers are  equal to $(-1)^{k-1}C_{k-1}$ with the Catalan numbers $$C_n :=\frac{\binom{2 n}{n}}{n+1}$$ 
 which are integers. The reduction to $a=\frac 12$ is done as for the Pascal triangle using the relation 
 $$
 \binom{x}{k}=\binom{x-1}{k}+\binom{x-1}{k-1}
 $$
 It is important to note, in order to perform the induction, that all the binomial coefficients $\binom{x}{k}$ with $x$ a half integer and $k<0$ are equal to zero.\newline
$(ii)$~One has 
$$
\binom{2\ell}{\ell}\, p_\ell(m)=\sum_{k=0}^\ell  2^{3 k}\binom{2\ell}{\ell}\,\binom{\ell}{k}\binom{2k}{k}^{-1}\binom{m}{k}
$$
By $(i)$ each term $\binom{2\ell}{\ell}\,\binom{\ell}{k}\binom{2k}{k}^{-1}$ is an integer thus since the other terms of the product are integers so are each term of the sum.\newline
$(iii)$~One has by \eqref{pmzg2} that $p_\ell(m)=p_m(\ell)$, thus one gets, using $(ii)$,
$$
\sigma(\ell,m)=\binom{2\ell}{\ell}\, p_\ell(m)\binom{2m}{m}\, p_m(\ell)\in \Z
$$
which gives the required integrality. \endproof 
 The reason why the expression $\sigma(\ell,m)$ captures the key structure appearing when one adjoins a prime to the Archimedean place is the following.
 \begin{prop}\label{matrixentries} The matrix coefficients $\gamma_{i,j}(q)$ of the matrix $\cA_n^+(q)\oplus \cA_n^-(q)$ are zero when $i-j$ is odd and are otherwise given by the Lambert series 
 	\begin{equation}\label{gammaij}
\gamma_{i,j}(q)=\sum_{\ell=0}^\infty (-1)^{j+\ell} 2^{\frac 32-2\ell-i-j} \sqrt{\sigma(\ell,i)\sigma(\ell,j)}	\frac{q^{2\ell+1}}{1-q^{2\ell+1}}
 \end{equation}
 \end{prop}
\proof One has by \eqref{pmzg}, the components $\eta_\ell(m)$ of the vectors $\eta_\ell$ involved in $T_\ell^\pm:=\pm 2^{3/2}\vert \eta_\ell^\pm\rangle \langle \eta_\ell^\pm\vert$ to obtain the Lambert series \eqref{lambertrank1} of rank one operators. By \eqref{pmzg1} and \eqref{pmzg2} one has
\begin{equation}\label{gammaij1}
\eta_\ell(m)=\frac{2^{-m}\sqrt{(2 m)!}}{m!}p_\ell(m)=2^{-m}\sqrt{\binom{2m}{m}}p_\ell(m)
 \end{equation}
 This gives 
 $$
 \alpha_\ell\, \eta_\ell(i)\,\eta_\ell(j)=(-1)^{\ell}2^{-2\ell}\binom{2\ell}{\ell}\times 2^{-i}\sqrt{\binom{2i}{i}}p_\ell(i)\times 2^{-j}\sqrt{\binom{2j}{j}}p_\ell(j)=$$ $$
 (-1)^{\ell}2^{-2\ell-i-j}\sqrt{\sigma(\ell,i)\sigma(\ell,j)}
 $$
 Finally the sign $(-1)^j$ and the overall factor $2^{3/2}$ account for \eqref{tlplusminus}.\endproof 
 The fundamental fact now is that the terms $\sqrt{\sigma(\ell,i)\sigma(\ell,j)}$ are cyclotomic integers, \ie 
 \begin{equation}\label{gammaij2}
 \sqrt{\sigma(\ell,i)\sigma(\ell,j)}\in \ccyc
  \end{equation}
  The next result gives in particular a proof of Fact \ref{keyfact}.
  \begin{thm}\label{integralitymain} $(i)$~The coefficients of the $q$-expansion of the  $\gamma_{i,j}(q)$ belong to the ring $\ccyc[\frac 12]$. \newline 
  $(ii)$~The coefficients of the $q$-expansion of $\Tr(\wedge^m \cA_n^\pm)$  belong to the ring $\Z[1/\sqrt 2]$.\newline
  $(iii)$~The coefficients of the $q$-expansion of $D_n(q)/D_n(0)$  belong to the ring $\Z[1/\sqrt 2]$.  	
  \end{thm} 
  \proof $(i)$~The coefficient $c_n(i,j)$ of $q^n$ in the $q$-expansion of  $\gamma_{i,j}(q)$ is obtained as the sum  $$\sum_{\ell, 2\ell+1\vert n}(-1)^{j+\ell} 2^{\frac 32-2\ell-i-j} \sqrt{\sigma(\ell,i)\sigma(\ell,j)}	$$ over the integers $\ell$ such that $2\ell +1$ divides $n$. Each term $\sqrt{\sigma(\ell,i)\sigma(\ell,j)}$ is a cyclotomic integer and thus $c_n(i,j)\in \ccyc[\frac 12]$.\newline
  $(ii)$~By Proposition \ref{similar}, the matrix $ \cA_n^\pm$ is similar to the restriction to even and odd of the matrix $\cD_n(0)^{-1}\cE_n(q)$ which has coefficients in $\Q[\sqrt 2,q]$, thus the coefficients of the $q$-expansion of $\Tr(\wedge^m \cA_n^\pm)$  belong to $\Q[\sqrt 2]$. The intersection of the rings $\Q(\sqrt 2)$ and $\ccyc[\frac 12]$ is the ring $\Z[1/\sqrt 2]$.\newline
  $(iii)$~Use $(ii)$ and \eqref{Ln0}. \endproof 
  The cyclotomic integers $\sqrt{\sigma(\ell,i)\sigma(\ell,j)}$ play a key role, and one can use  \eqref{halfinteg} to obtain 
  $$
  \binom{2\ell}{\ell}\, p_\ell(m)=\sum_{k=0}^\ell  2^{ k+2\ell}\binom{\ell-\frac 12}{\ell-k}\binom{m}{k}
$$
so that 
\begin{equation}\label{halfinteg1}
\sigma(a,b)=\sum_{k=0}^a\sum_{k'=0}^b  2^{k+ k'+2a+2b}\binom{a-\frac 12}{a-k}\binom{b-\frac 12}{b-k'}\binom{b}{k}\binom{a}{k'}
 \end{equation}
\section{Eigenvalues of the $q$-matrices $\cA_n(q)$} \label{eigenv}
When one computes the trace of the exterior powers $\Tr(\wedge^m \cA_n^\pm)$, one observes
\begin{fact}\label{keyfact2} The $q$-expansion of the trace $\Tr(\wedge^m \cA_n^\pm)$
starts by the power $q^{m^2}$.
\end{fact} 
This fact suggests to look for $q$-expansions of eigenvalues $\lambda_n(k)$, $0\leq k\leq n$, of the form
\begin{equation}\label{ansatzeigen}
\lambda_n^\pm(k)=\sum_{j\geq 2k+1} \lambda_n^\pm(k,j)q^j
 \end{equation}
When this is done, one finds  for $n=1$ the two eigenvalues $\sqrt 2\times \lambda_1^+(k)$, $k\in \{0,1\}$ for $\cA_1^+$
$$
\frac{11}{2}q+\frac{11}{2}q^2-\frac{983}{44}q^3+\frac{11}{2}q^4+\frac{6886331}{21296}q^5-\frac{1059743}{5324}q^6-\frac{15326598829}{5153632}q^7+\frac{1001841961}{322102}q^8+\ldots
$$
and 
$$
-\frac{384}{11}q^3+\frac{152000}{1331}q^5+\frac{188736}{1331}q^6+\frac{234365040}{161051}q^7-\frac{500035200}{161051}q^8-\frac{281727520584}{19487171}q^9+\ldots
$$
In this case $n=1$ one has two eigenvalues $\lambda_1^+(k)$, $k\in \{0,1\}$ and one determines the coefficients $\lambda_1^+(k,j)$ using the coefficients in
\begin{equation}\label{tracewedge}
\Tr(\wedge^m \cA_n^+)=\sum_{j\geq m^2} \gamma_n(m,j)q^j
\end{equation}
and the  equalities coming from $\Tr(\wedge^1 \cA_1^+)=\lambda_1^+(0)+\lambda_1^+(1)$, $\Tr(\wedge^2 \cA_1^+)=\lambda_1^+(0)
\lambda_1^+(1)$
\begin{equation}\label{tracewedge1}
\lambda_1^+(0,j)+\lambda_1^+(1,j)=\gamma_1(1,j), \ \ \sum_{i+i'=j} \lambda_1^+(0,i)\lambda_1^+(1,i')=\gamma_1(2,j)
\end{equation}
One has $\lambda_1^+(1,j)=0$ for $j<3$ thus \eqref{tracewedge1} determines $\lambda_1^+(0,j)$ for $j<3$. Note that $\lambda_1^+(0,0)=0$ and that $\lambda_1^+(0,1)\neq 0$. One then uses \eqref{tracewedge1} to get $\lambda_1^+(0,1)\lambda_1^+(1,3)=\gamma_1(2,4)$ which determines $\lambda_1^+(1,3)$. Thus one gets $\lambda_1^+(0,3)$ since $\lambda_1^+(0,3)+\lambda_1^+(1,3)=\gamma_1(1,3)$. Then the equality 
$$
\lambda_1^+(0,1)\lambda_1^+(1,4)+\lambda_1^+(0,2)\lambda_1^+(1,3)=\gamma_1(2,5)
$$
determines $\lambda_1^+(1,4)$ and hence $\lambda_1^+(0,4)$. By induction one then determines $\lambda_1^+(1,j)$ using
$$
\lambda_1^+(0,1)\lambda_1^+(1,j)+\sum_{i+i'=j+1,i>1} \lambda_1^+(0,i)\lambda_1^+(1,i')=\gamma_1(2,j+1)
$$
This, in turns, determines the $\lambda_1^+(0,j)$ by $\lambda_1^+(0,j)+\lambda_1^+(1,j)=\gamma_1(1,j)$.\newline
For $n=2$ one has three eigenvalues $\lambda_2^+(k)$, $k\in \{0,1,2\}$. Since $\lambda_2^+(2)$ starts by the monomial $\lambda_2^+(2,5)q^5$ it does not interfere with the computation as above of the coefficients $\lambda_2^+(k,j)$ for $j<5$. One gets 
\begin{equation}\label{tracewedge2}
\lambda^+ _2(0,5)+\lambda^+ _2(1,5)+\lambda^+ _2(2,5)= \gamma_2(1,5), \  \ \lambda^+ _2(0,2) \lambda^+ _2(1,3)+\lambda^+ _2(0,1) \lambda^+ _2(1,4)= \gamma_2(2,5)
\end{equation}
and 
\begin{equation}\label{tracewedge3}
\lambda^+ _2(0,3) \lambda^+ _2(1,3)+\lambda^+ _2(0,2) \lambda^+ _2(1,4)+\lambda^+ _2(0,1) \lambda^+ _2(1,5)+\lambda^+ _2(0,1) \lambda^+ _2(2,5)=\gamma_2(2,6)
\end{equation}
\begin{equation}\label{tracewedge4}
\lambda^+ _2(0,1) \lambda^+ _2(1,3) \lambda^+ _2(2,5)=\gamma_2(3,9)
\end{equation}
Thus \eqref{tracewedge4} determines $\lambda^+ _2(2,5)$, then \eqref{tracewedge3} determines $\lambda^+ _2(1,5)$ and finally \eqref{tracewedge2} determines $\lambda^+ _2(0,5)$. The next step uses 
$$
\lambda^+ _2(0,2) \lambda^+ _2(1,3) \lambda^+ _2(2,5)+\lambda^+ _2(0,1) \lambda^+ _2(1,4) \lambda^+ _2(2,5)+\lambda^+ _2(0,1) \lambda^+ _2(1,3) \lambda^+ _2(2,6)=\gamma_2(3,10)
$$
which determines $ \lambda^+ _2(2,6)$. One then uses 
$$
\lambda^+ _2(0,4) \lambda^+ _2(1,3)+\lambda^+ _2(0,3) \lambda^+ _2(1,4)+\lambda^+ _2(0,2) \lambda^+ _2(1,5)+$$ $$+\lambda^+ _2(0,1) \lambda^+ _2(1,6)+\lambda^+ _2(0,2) \lambda^+ _2(2,5)+\lambda^+ _2(0,1) \lambda^+ _2(2,6)=\gamma_2(2,7)
$$
to determine $\lambda^+ _2(1,6)$, and $\lambda^+ _2(0,6)+\lambda^+ _2(1,6)+\lambda^+ _2(2,6)=\gamma_2(1,6)$ to get $\lambda^+ _2(0,6)$.\newline 
The next step 
$$
\lambda^+ _2(0,3) \lambda^+ _2(1,3) \lambda^+ _2(2,5)+\lambda^+ _2(0,2) \lambda^+ _2(1,4) \lambda^+ _2(2,5)+\lambda^+ _2(0,1) \lambda^+ _2(1,5) \lambda^+ _2(2,5)+$$ $$+\lambda^+ _2(0,2) \lambda^+ _2(1,3) \lambda^+ _2(2,6)+\lambda^+ _2(0,1) \lambda^+ _2(1,4) \lambda^+ _2(2,6)+\lambda^+ _2(0,1) \lambda^+ _2(1,3) \lambda^+ _2(2,7)=\gamma_2(3,11)
$$
determines $\lambda^+ _2(2,7)$ and one then uses 
$$
\lambda^+ _2(0,5) \lambda^+ _2(1,3)+\lambda^+ _2(2,5) \lambda^+ _2(1,3)+\lambda^+ _2(0,4) \lambda^+ _2(1,4)+\lambda^+ _2(0,3) \lambda^+ _2(1,5)+\lambda^+ _2(0,2) \lambda^+ _2(1,6)+$$ $$+\lambda^+ _2(0,1) \lambda^+ _2(1,7)+\lambda^+ _2(0,3) \lambda^+ _2(2,5)+\lambda^+ _2(0,2) \lambda^+ _2(2,6)+\lambda^+ _2(0,1) \lambda^+ _2(2,7)=\gamma_2(2,8)
$$
to get $\lambda^+ _2(1,7)$, and finally $\lambda^+ _2(0,7)$ using $\gamma_2(1,7)$. 

\begin{lem}\label{usefact2} Assume Fact \ref{keyfact2}. Then the eigenvalues of the matrix $\cA_n^\pm(q)$ are given by power series of the form \eqref{ansatzeigen} where the coefficients $\lambda^\pm_n(k,j)$ are rational number multiples of $\sqrt 2$.	
\end{lem}
\proof We first define for each $n$ a sequence $s_n(r)$, $r\in \N$, where, for each $r$, $s_n(r)=(k,j)$ is a pair of positive integers where the index $k\geq 0$ indicates which eigenvalue is considered, and $j\geq 2k+1$ indicates the power of $q$ involved. The sequence $s_n(r)$ governs the inductive computation of the coefficients $\lambda_n(s_n(r))$. The computation of the next one $\lambda _n(s_n(r+1))$ takes into account that all previous $\lambda _n(s_n(r'))$, $r'\leq r$ have been determined. The index $j$ in $s_n(r)=(k,j)$ is non-decreasing in $r$, but for fixed $j$ the index $k$ is strictly decreasing so that if $k>0$ one has 
$$
s_n(r)=(k,j)\Rightarrow s_n(r+1)=(k-1,j)
$$
Thus, to determine the sequence $s_n(r)$ one needs to know what is $s_n(r+1)$ when $s_n(r)=(0,j)$. This depends on $j$ as follows, where $E(x)$ is the integer part of $x$, 
$$
s_n(r)=(0,j)\Rightarrow s_n(r+1)=(\min\{n,(E(j/2)\},j+1)
$$
We then define another sequence $t_n(r)$, $r\in \N$, where, for each $r$, $t_n(r)=(k,j)$ is a pair of positive integers where the index $k\geq 0$ indicates which symmetric function of the eigenvalues is considered, and $j$ indicates the power of $q$ involved. The first term $k$ in $t_n(r)$ is $1$ plus the term $k$ in $s_n(r)$ because the indexation of symmetric functions differs by $1$ of the indexation of the eigenvalues. To obtain the second term $j$ in $t_n(r)=(k,j)$ one needs to minimise the degree of a cofactor of $\lambda _n(s_n(r))$ in a symmetric function involving a product of $k$ terms. Such a cofactor is a product of $k-1$ terms of the form $\prod \lambda _n(\ell,i)$ and the $\ell$ are distinct while their degrees
are as small as possible \ie are $2\ell+1$. Thus the degree of the cofactor is 
$$
\sum_0^{k-1} (2\ell+1)=k^2
$$ 
and one has 
$$
s_n(r)=(k,j)\Rightarrow t_n(r)=(k+1,k^2+j)
$$
The inductive process is done for increasing $r$. One determines inductively, for fixed $n$, the value of $\lambda_n(s_n(r))$ using the equality of the product of $\lambda_n(s_n(r))=\lambda_n(k,j)$ by its cofactor in the $k+1$ symmetric function added to the other terms in this symmetric function, with the value of the symmetric function computed as the trace $\Tr(\wedge^{k+1} \cA_n^\pm)$. The exponent of $q$ that one considers in the equality is the second term in $t_n(r)=(k+1,k^2+j)$ and is thus in general far larger than $j$. Thus computing the term of exponent $j$, \ie $\lambda_n(k,j)$, requires computing the terms of higher exponents in $\Tr(\wedge^{k+1} \cA_n^\pm)$. When dealing with these terms one needs to make sure that except for the term involving $\lambda_n(k,j)$, all other terms are already computed.
\subsection{Proof of Fact \ref{keyfact2}} 
We first prove a general Lemma which describes the trace $\tr\left(\wedge^m(T)\right)$ for an operator $T$ given as a sum of rank one operators.
\begin{lem}\label{wedgetrace} Let $\cH$ be a  Hilbert space and $\iota:\cH\otimes \bar \cH\to \End(\cH)$ the linear map 
$$
\iota(\xi\otimes \bar \eta)(\alpha):=\xi \langle \eta,\alpha\rangle, \ \ \forall \alpha \in \cH
$$
$(i)$~The following equality defines a Hilbert space structure on $\wedge^m(\cH)$ viewed as a quotient of $\otimes^m\cH$,
	$$
	\langle \xi_1\otimes \xi_2\otimes \ldots \otimes \xi_m, \eta_1\otimes \eta_2\otimes \ldots \otimes \eta_m\rangle=\det \left(\langle \xi_i,\eta_j\rangle\right)
	$$
	$(ii)$~Let $T=\sum_{k\in I} \iota(\xi_k\otimes \bar \eta_k)\in \End(\cH) $. One has 
	$$
	\wedge^m(T)=\frac{1}{m!}\sum \iota\left(\xi_{k_1}\otimes \xi_{k_2}\otimes \ldots \otimes \xi_{k_m}\otimes \bar\eta_{k_1}\otimes \bar\eta_{k_2}\otimes \ldots \otimes \bar\eta_{k_m} \right) 
	$$
	where the sum is over all choices of maps $i\mapsto k_i$ from $\{1,\ldots,m\}$ to the index set $I$.\newline
	$(iii)$~Let $T=\sum_{k\in I} \iota(\xi_k\otimes \bar \eta_k)\in \End(\cH) $. Then 
	$$
	\tr\left(\wedge^m(T)\right)=\frac{1}{m!}\sum \det \left(\langle \xi_{k_i},\eta_{k_j}\rangle\right)
	$$
	\end{lem}
\proof $(i)$~Let $(e_j)_{j\in \N}$ be an orthonormal basis of $\cH$, then the vectors of the form 
$$
e(J):=e_{i_1}\otimes \ldots \otimes e_{i_m}, \  \ i_j\in J, \  i_1<i_2<\ldots <i_m
$$
form an orthonormal basis of $\wedge^m(\cH)$.\newline
$(ii)$~The normalization factor $\frac{1}{m!}$ comes  from the equality of operators in 
$\wedge^m(\cH)$,
$$
\iota\left(\xi_{1}\otimes \ldots \otimes \xi_{m}\otimes \bar\eta_{1}\otimes  \ldots \otimes \bar\eta_{m} \right) =\iota\left(\xi_{k_1}\otimes \ldots \otimes \xi_{k_m}\otimes \bar\eta_{k_1}\otimes  \ldots \otimes \bar\eta_{k_m} \right) 
$$
for any permutation $k$ of $\{1,\ldots,m\}$.\newline
$(iii)$~Follows from $(ii)$ and $(i)$. \endproof 
Let us now prove Fact \ref{keyfact2}.
\begin{lem}
\label{keylem2} The $q$-expansion of the trace $\Tr(\wedge^m \cA_n^\pm)$
starts by  $c_{n,m}^\pm \,q^{m^2}$ with coefficient given by 
\begin{equation}\label{coeffqm2}
c_{n,m}^+=(-1)^{\frac{m(m-1)}{2}}2^{3m/2}\ \det_{0\leq \ell,\ell'<m} \left(\langle \xi_{n,\ell}^+\vert \xi_{n,\ell'}^+\rangle\right)
\end{equation}
where the vectors $\xi_{n,\ell}^+$ have the coordinates $2^{-\ell-2j}\sqrt{\sigma(\ell,2j)}$ for $0\leq 2j\leq n$. In the odd case one gets
\begin{equation}\label{coeffqm2odd}
c_{n,m}^-=(-1)^{\frac{m(m+1)}{2}}2^{3m/2}\ \det_{0\leq \ell,\ell'<m} \left(\langle \xi_{n,\ell}^-\vert \xi_{n,\ell'}^-\rangle\right)
\end{equation}
where the vectors $\xi_{n,\ell}^-$ have the coordinates $2^{-\ell-2j-1}\sqrt{\sigma(\ell,2j+1)}$ for $0\leq 2j+1\leq n$.
\end{lem} 
\proof By Proposition \ref{matrixentries}, the matrix $\cA_n^+$ has entries $\gamma_{2i,2j}^+(q)$ given by
$$
\gamma_{2i,2j}^+(q)=\sum_{\ell=0}^\infty (-1)^{\ell} 2^{\frac 32-2\ell-2i-2j} \sqrt{\sigma(\ell,2i)\sigma(\ell,2j)}	\frac{q^{2\ell+1}}{1-q^{2\ell+1}}
$$
It is the sum of rank one matrices $$\sum_{\ell=0}^\infty (-1)^{\ell} 2^{\frac 32}\frac{q^{2\ell+1}}{1-q^{2\ell+1}}\ \xi_{n,\ell}^+\otimes \xi_{n,\ell}^+$$
where the vectors $\xi_{n,\ell}^+$ have the coordinates $2^{-\ell-2j}\sqrt{\sigma(\ell,2j)}$ for $0\leq 2j\leq n$. The smallest power of $q$ which appears in $\Tr(\wedge^m \cA_n^\pm)$ comes, by Lemma \ref{wedgetrace} $(iii)$, from the  choice of $m$ different $\ell_i$ with smallest $\sum (2\ell_i+1)$. This choice is $\ell_1=0$, $\ell_i=i-1$, $\ell_m=m-1$. It gives $\sum (2\ell_i+1)=m^2$, the product of the $(-1)^{\ell_i}$ is $(-1)^{\frac{m(m-1)}{2}}$, thus one gets \eqref{coeffqm2}. In the odd case one has 
$$
\gamma_{2i+1,2j+1}^-(q)=-\sum_{\ell=0}^\infty (-1)^{\ell} 2^{\frac 32-2\ell-2i-1-2j-1} \sqrt{\sigma(\ell,2i+1)\sigma(\ell,2j+1)}	\frac{q^{2\ell+1}}{1-q^{2\ell+1}}
$$
It is the sum of rank one matrices $$\sum_{\ell=0}^\infty (-1)^{\ell+1} 2^{\frac 32}\frac{q^{2\ell+1}}{1-q^{2\ell+1}}\ \xi_{n,\ell}^-\otimes \xi_{n,\ell}^-$$
where the vectors $\xi_{n,\ell}^-$ have the coordinates $2^{-\ell-2j-1}\sqrt{\sigma(\ell,2j+1)}$ for $0\leq 2j+1\leq n$.\endproof 
The first values of the coefficients $c^+_{n,m}$ are the following 
$$
\left(
\begin{array}{cccc}
 2 \sqrt{2} & 0 & 0 & 0 \\
 \frac{11}{2 \sqrt{2}} & -96 & 0 & 0 \\
 \frac{211}{32 \sqrt{2}} & -\frac{1609}{4} & -35840 \sqrt{2} & 0 \\
 \frac{1919}{256 \sqrt{2}} & -\frac{524501}{512} & -504077 \sqrt{2} & 423886848 \\
\end{array}
\right)
$$
The first values of the coefficients $c^-_{n,m}$ are the following
$$
\left(
\begin{array}{cccc}
 -\sqrt{2} & 0 & 0 & 0 \\
 -\frac{13}{4 \sqrt{2}} & -40 & 0 & 0 \\
 -\frac{271}{64 \sqrt{2}} & -\frac{2971}{16} & 13440 \sqrt{2} & 0 \\
 -\frac{2597}{512 \sqrt{2}} & -\frac{1038803}{2048} & \frac{1624893}{4 \sqrt{2}} & 147603456 \\
\end{array}
\right)
$$

\subsection{Integrality of eigenvalues} We can now combine Lemma \ref{keylem2} with Lemma \ref{usefact2} to obtain
\begin{thm}\label{eigenthm}
The eigenvalues of the matrix $\cA_n^\pm(q)$ are given by power series of the form \eqref{ansatzeigen} where the coefficients $\lambda^\pm_n(k,j)$ are rational number multiples of $\sqrt 2$.	
\end{thm}
The first examples of $q$-expansions of eigenvalues are, first for $n=0$ the eigenvalue
$$
\lambda^+_0(0)=2 \sqrt{2} q+2 \sqrt{2} q^2+\sqrt{2} q^3+2 \sqrt{2} q^4+\frac{11 q^5}{2 \sqrt{2}}+\sqrt{2} q^6+\ldots 
$$
Next, for $n=1$ the two eigenvalues
$$
\lambda^+_1(0)=\frac{11 q}{2 \sqrt{2}}+\frac{11 q^2}{2 \sqrt{2}}-\frac{983 q^3}{44
   \sqrt{2}}+\frac{11 q^4}{2 \sqrt{2}}+\frac{6886331 q^5}{21296 \sqrt{2}}+ \ldots,$$ $$\lambda^+_1(1)=-\frac{192}{11}
   \sqrt{2} q^3+\frac{76000 \sqrt{2} q^5}{1331}+\ldots
   $$
Next for $n=2$, the three eigenvalues
$$
\lambda^+_2(0)=\frac{211 q}{32 \sqrt{2}}+\frac{211 q^2}{32 \sqrt{2}}-\frac{1244983 q^3}{13504
   \sqrt{2}}+\frac{211 q^4}{32 \sqrt{2}}+\frac{7554810186491 q^5}{2404846336
   \sqrt{2}}-\frac{1154408355343 q^6}{601211584 \sqrt{2}}\ldots$$ $$
   \lambda^+_2(1)=-\frac{12872}{211} 
   \sqrt{2} q^3+\frac{7392097696500 \sqrt{2} q^5}{15114834979}+\frac{8012710588 \sqrt{2}
   q^6}{9393931}+\ldots$$ $$\lambda^+_2(2)=\frac{143360 \sqrt{2} q^5}{1609}-\frac{482941605888 \sqrt{2}
   q^7}{4165509529}-\frac{20468604828168211968 \sqrt{2} q^9}{10784008474947049}+\ldots
   $$
Note that the denominators are no longer powers of $2$ and for  instance the last denominator above is the fifth power of the prime $1609$. For $n=3$ one gets
$$
\lambda^+_3(0)=\frac{1919 q}{256 \sqrt{2}}+\frac{1919 q^2}{256 \sqrt{2}}-\frac{218085103 q^3}{982528
   \sqrt{2}}+\frac{1919 q^4}{256 \sqrt{2}}+\frac{204982567422146795 q^5}{14472877176832
   \sqrt{2}}$$ $$-\frac{31074822021446383 q^6}{3618219294208
   \sqrt{2}}-\frac{74352350022892667863024709 q^7}{106594506098383253504 \sqrt{2}}+\ldots
   $$
   $$
  \lambda^+_3(1)= -\frac{524501 q^3}{1919 \sqrt{2}}+\frac{29149147149228946125 q^5}{7413123586060118
   \sqrt{2}}+\frac{114385854608603 q^6}{14133669118
   \sqrt{2}}+$$ $$+\frac{10190268835092002960940341810228356358125
   q^7}{30040269345915199837200814357616792 \sqrt{2}}+\ldots
   $$
   $$
  \lambda^+_3(2)= \frac{258087424 \sqrt{2} q^5}{524501}-\frac{13031546695237024906049280 \sqrt{2}
   q^7}{10390532462737992630511}$$ $$-\frac{511662216808719861706456813376030807133551384256
   \sqrt{2} q^9}{14822727630026862171674605299320669826184831}+\ldots
   $$
   $$
  \lambda^+_3(3)= -\frac{30277632 \sqrt{2} q^7}{72011}+\frac{116176571793604608 \sqrt{2}
   q^9}{373419098137331}+\ldots
   $$
   For the eigenvalues of the matrix $\cA_n^-(q)$ one gets the following expressions
   $$
   \lambda^-_0(0)=-\sqrt{2} q-\sqrt{2} q^2+\frac{23 q^3}{\sqrt{2}}-\sqrt{2} q^4-\frac{251 q^5}{4 \sqrt{2}}+\frac{23
   q^6}{\sqrt{2}}+\ldots
   $$ 
   Next, for $n=1$ the two eigenvalues
$$
\lambda^-_1(0)=-\frac{13 q}{4 \sqrt{2}}-\frac{13 q^2}{4 \sqrt{2}}+\frac{10687 q^3}{104 \sqrt{2}}-\frac{13 q^4}{4
   \sqrt{2}}-\frac{121293131 q^5}{70304 \sqrt{2}}+\frac{15918103 q^6}{17576 \sqrt{2}},$$ $$\lambda^-_1(1)=\frac{160 \sqrt{2}
   q^3}{13}+\frac{151920 \sqrt{2} q^5}{2197}-\frac{854960 \sqrt{2} q^6}{2197}   $$
   Next for $n=2$, the three eigenvalues
   $$
   \lambda^-_2(0)=-\frac{271 q}{64 \sqrt{2}}-\frac{271 q^2}{64 \sqrt{2}}+\frac{8871127 q^3}{34688 \sqrt{2}}-\frac{271 q^4}{64
   \sqrt{2}}-\frac{104586068271995 q^5}{10190085632 \sqrt{2}}+\frac{14369266864375 q^6}{2547521408
   \sqrt{2}},$$
   $$\lambda^-_2(1)=\frac{11884 \sqrt{2} q^3}{271}+\frac{22796026081686 \sqrt{2} q^5}{59130360181}-\frac{52712236634
   \sqrt{2} q^6}{19902511}+\ldots,$$ $$\lambda^-_2(2)=-\frac{215040 \sqrt{2} q^5}{2971}+\ldots
   $$
   One notices that in all these expressions for the eigenvalue $\lambda^\pm_n(0)$, the coefficients of $q$ and $q^2$ are the same, as follows from the basic property of Lambert series involving only odd powers. Moreover by \eqref{an2bis} the determinant $\Delta_n(q)=\sum_0^\infty \mu_{n,k}\ q^k$ is the product 
   \begin{equation}\label{plusminus}
   \Delta_n(q)=\prod_{j=0}^a (1+\lambda^+_a(j))\prod_{k=0}^b (1+\lambda^-_b(k))
   \end{equation}
   where the integers $a,b$ are determined as the dimensions of the subspaces $\cH_n^\pm$ by $1+a=\dim( \cH_n^+)$ and $1+b=\dim( \cH_n^-)$, \ie $a+1$ is the number of even elements in the set $\{0, 1, \ldots,n\}$ and $b+1$ the number of odd elemnts. Thus in particular one has $a+b=n-1$.
   \begin{prop} \label{coeffqeigen} $(i)$~The coefficient $\gamma_n^\pm$  of $q$ in the eigenvalue $\lambda^\pm_n(0)$ is equal to the coefficient of $q^2$ and is given by 
   	\begin{equation}\label{gammanpm}
  \gamma_n^+=2 \sqrt{2}\sum_{j=0}^n \alpha_{2j}, \  \  \gamma_n^-=2 \sqrt{2}\sum_{j=0}^n \alpha_{2j+1}
   \end{equation}
   $(ii)$~The coefficient $\mu_{n,2}$ of $q^2$ in $\Delta_n(q)$ is given, according to the parity of $n$, by 
   \begin{equation}\label{coeffq2}
  \mu_{2m+1,2}=\gamma_m^++\gamma_m^-+\gamma_m^+\gamma_m^-, \  \  \mu_{2m,2}=\gamma_m^++\gamma_{m-1}^-+\gamma_m^+\gamma_{m-1}^-
   \end{equation}
   \end{prop}
   \proof $(i)$~Follows from Proposition \ref{similar}.\newline
   $(ii)$~Follows from \eqref{plusminus}.\endproof

\subsection{Orthogonal polynomials $P_n(x,q)$}
We are now ready to prove the main result
\begin{thm} \label{qthm} Let $p$ be a prime, and $d\mu$ be the measure on $\R$ given by
$$
d\mu(s)=\vert\frac{1}{1- p^{-\frac{1}{2}-is}}\Gamma(\frac 14+\frac 12 is)\vert^2 ds
$$
The coefficients of the characteristic polynomial of the compression of the operator of multiplication by $x$ on the subspace of $L^2(\R, d\mu)$ given by polynomials of degree $\leq n$ are universal $q$-series with coefficients in the ring $\Z[1/\sqrt 2]$ evaluated at $q=1/p$.	
\end{thm}
\proof This follows from the same property for the coefficients $a_n^2$ and the formula of section \ref{sectorthopoly} for the Jacobi matrix in terms of the $a_n^2$.\endproof 
We display a few of these polynomials, one has $P_1(x,q)=x$, 
$$
P_2(x,q)=-\frac{1}{2}+x^2+\frac{3 q}{\sqrt{2}}+\left(\frac{3}{\sqrt{2}}-6\right) q^2+\left(\frac{27}{2 \sqrt{2}}-12\right) q^3+$$ $$+\left(\frac{75}{\sqrt{2}}-36\right) q^4+\left(\frac{1905}{8 \sqrt{2}}-186\right) q^5+\left(\frac{1851}{2 \sqrt{2}}-654\right) q^6+\ldots
$$
$$
P_3(x,q)=-\frac{7}{2} x+x^3+x\biggl(-\frac{15 q }{2 \sqrt{2}}+ \left(-\frac{15 }{2}-\frac{15 }{2 \sqrt{2}}\right)q^2+ \left(\frac{885 }{4 \sqrt{2}}-15 \right)q^3+$$ $$+ \left(300 -\frac{105 }{2 \sqrt{2}}\right)q^4+ \left(\frac{495 }{2}-\frac{11325 }{16 \sqrt{2}}\right)q^5+ \left(\frac{6645 }{4 \sqrt{2}}-\frac{7095 }{2}\right)q^6+\ldots\biggl)
$$

\section{Diagonal matrix $D$ and hypergeometric functions} \label{diag}
We recall that in the context of general coefficients $a(n)$ for the Jacobi matrix $A$, we showed in \cite{c2m}, Lemma 5.1, that one can use the integer grading $N$ to define
\begin{equation}\label{eplusminusguess}
E_+:=\frac 14+N+\frac{1}{2i}[A,N],\ \ E_-:=\frac 14+N-\frac{1}{2i}[A,N]		
\end{equation}
One then has (taking into account the added $\frac 14$)
\begin{enumerate}
	\item The commutator $[E_+,E_-]$ is equal to $-iA$.
	\item The commutator $[A,E_+]$ is equal to $2i E_++D$ where the diagonal matrix $D$ has the entries $d_n$ with 
\begin{equation}\label{dn8}
	-i\ d_n=-2n-\frac 12+a(n)^2-a(n-1)^2
\end{equation}
\end{enumerate}
 We compute the first examples of the $q$-expansion of the $-i \ d_n$, we know that they are $q$-series with coefficients in the ring $\Z[1/\sqrt 2]$,
 $$
-i\ d_1= \frac{27 q}{2 \sqrt{2}}-\frac{9 q^2}{2}+\frac{27 q^2}{2 \sqrt{2}}-9 q^3-\frac{777 q^3}{4 \sqrt{2}}-372
   q^4+\frac{405 q^4}{2 \sqrt{2}}$$ $$-\frac{1239 q^5}{2}+\frac{18945 q^5}{16 \sqrt{2}}+\frac{4479 q^6}{2}+\frac{759
   q^6}{4 \sqrt{2}}+\ldots
   $$
 $$
-i\ d_2= -\frac{249 q}{8 \sqrt{2}}+\frac{867 q^2}{32}-\frac{249 q^2}{8 \sqrt{2}}+\frac{867 q^3}{16}+\frac{78171 q^3}{64
   \sqrt{2}}-\frac{51909 q^4}{256}-\frac{50475 q^4}{64 \sqrt{2}}+$$ $$+\frac{97029 q^5}{64}-\frac{5044725 q^5}{512
   \sqrt{2}}+\frac{47751987 q^6}{2048}+\frac{2469099 q^6}{512 \sqrt{2}}+\ldots
   $$
 $$
 -i\ d_3=\frac{855 q}{16 \sqrt{2}}-\frac{4185 q^2}{128}+\frac{855 q^2}{16 \sqrt{2}}-\frac{4185 q^3}{64}-\frac{2237325
   q^3}{512 \sqrt{2}}-\frac{11265225 q^4}{4096}+$$ $$+\frac{819945 q^4}{512 \sqrt{2}}-\frac{5638785
   q^5}{1024}+\frac{1028457435 q^5}{16384 \sqrt{2}}+\frac{19805812335 q^6}{131072}-\frac{667328265 q^6}{16384
   \sqrt{2}}+\ldots
 $$
 \subsection{Coefficient of $q$ in $-i\ d_n$}
The list of coefficients of $q$ in the first $-i\ d_n$ is as follows
$$
\frac{27}{2 \sqrt{2}},-\frac{249}{8 \sqrt{2}},\frac{855}{16 \sqrt{2}},-\frac{10185}{128 \sqrt{2}},\frac{27909}{256 \sqrt{2}},-\frac{144837}{1024 \sqrt{2}},\frac{361647}{2048 \sqrt{2}},-\frac{7020585}{32768 \sqrt{2}},$$ $$\frac{16664505}{65536 \sqrt{2}},-\frac{77736087}{262144 \sqrt{2}},\frac{178738833}{524288 \sqrt{2}},-\frac{1624521717}{4194304 \sqrt{2}},\frac{3654510825}{8388608 \sqrt{2}},$$ $$-\frac{16302383325}{33554432 \sqrt{2}},\frac{36093907935}{67108864 \sqrt{2}},-\frac{1270383404265}{2147483648 \sqrt{2}}
$$
We now give their explicit general form
\begin{prop}\label{coeffq} The coefficient of $q$ in $-i\ d_n$ is given by 
	\begin{equation}\label{coeffqpert}
	(-1)^{n+1}\frac{ (16 n^2+8 n +3) \Gamma \left(n+\frac{1}{2}\right)}{\sqrt{2 \pi } \Gamma (n+1)}	=-\frac{1}{\sqrt 2}(16 n^2+8 n +3)\alpha_n
	\end{equation}
\end{prop}
\proof One has by \eqref{qexpandn1} and \eqref{mun1}, the equalities
\begin{equation}\label{qexpandn1bis}
\Delta_n(q)=D_n(q)/D_n(0)=1+\sum_1^\infty \mu_{n,k}\ q^k, \  \  \mu_{n,1}=2 \sqrt{2}\ \sum_0^n\alpha_k, \  \alpha_\ell=(- 4)^{-\ell}\binom{2\ell}{\ell}
\end{equation}
 Moreover by Corollary \ref{ansq1},
 $$
a(n)^2(q)=(n+\frac 12)(n+1)\frac{\Delta_{n-1}(q)\Delta_{n+1}(q)}{\Delta_n(q)^2},\ \ 
$$
which gives 
$$
a(n)^2(q)=(n+\frac 12)(n+1)\left(1+(\mu_{n-1,1}+\mu_{n+1,1}-2\mu_{n,1})q\right)+O(q^2)=
$$	
$$
=(n+\frac 12)(n+1)\left(1+2 \sqrt{2}\ (\alpha_{n+1}-\alpha_n)q\right)+O(q^2)
$$
so that,
$$
-2n-\frac 12+a(n)^2-a(n-1)^2=2 \sqrt{2}\left((n+\frac 12)(n+1)(\alpha_{n+1}-\alpha_n)- n(n-\frac 12)(\alpha_{n}-\alpha_{n-1})\right) 
$$
and  by \eqref{dn8},
$$
-i\ d_n=2 \sqrt{2}\left((n+\frac 12)(n+1)(\alpha_{n+1}-\alpha_n)- n(n-\frac 12)(\alpha_{n}-\alpha_{n-1})\right) 
$$
One has 
\begin{equation}\label{shiftan}	
(n-\frac 12)\,\alpha_{n-1}=-n\, \alpha_n, \ \ (n+1)\,\alpha_{n+1}=-(n+\frac 12)\,\alpha_{n}
\end{equation}
and thus the parenthesis can be written as a multiple of $\alpha_n$,
$$
\left((n+\frac 12)(n+1)(\alpha_{n+1}-\alpha_n)- n(n-\frac 12)(\alpha_{n}-\alpha_{n-1})\right)=p(n)\,\alpha_n
$$
$$
p(n)=-(n+\frac 12)^2-(n+\frac 12)(n+1)- n(n-\frac 12)-n^2=-\left(4n^2+2n+\frac 34 \right)
$$
which gives the result using $$\alpha_n=(- 4)^{-n}\binom{2n}{n}=\frac{(-1)^n \Gamma \left(n+\frac{1}{2}\right)}{\sqrt{\pi } \Gamma (n+1)}$$.\endproof
\subsection{Coefficient of $q^2$ in $-i\ d_n$} 
We now investigate the coefficient of $q^2$ in $-i\ d_n$. The main result is Theorem \ref{hinthyperg} below.  We use the notation of Proposition \ref{coeffqeigen}.
\begin{lem}\label{consratio} $(i)$~Let $n=2m$. The first terms of the $q$-expansion are  
\begin{equation}\label{qq3}	
	\Delta_{n-1}^+/\Delta_{n}^+=1-2 \sqrt{2}\,\alpha_{n}\ q-2 \sqrt{2}\,\alpha_{n}(1-\gamma^+_{m})\ q^2+0(q^3)
\end{equation}
\begin{equation}\label{qq3b}	
	\Delta_{n+1}^-/\Delta_{n}^-=1+2 \sqrt{2}\,\alpha_{n+1}\, q+2 \sqrt{2}\,\alpha_{n+1}(1-\gamma^-_{m-1})\ q^2+0(q^3)
\end{equation}
$(ii)$~Let $n=2m+1$. The first terms of the $q$-expansion are  
\begin{equation}\label{qq4}	
	 \Delta_{n+1}^+/\Delta_{n}^+=1+2 \sqrt{2}\,\alpha_{n+1}\ q+2 \sqrt{2}\,\alpha_{n+1}(1-\gamma^+_{m})\ q^2+0(q^3)
\end{equation}
\begin{equation}\label{qq4b}	
	\Delta_{n-1}^-/\Delta_{n}^-=1-2 \sqrt{2}\,\alpha_{n}\ q-2 \sqrt{2}\,\alpha_{n}(1-\gamma^-_{m})\ q^2+0(q^3)
\end{equation}
\end{lem}
\proof We  use Proposition \ref{coeffqeigen} which states that
the coefficient $\gamma_m^\pm$  of $q$ in the eigenvalue $\lambda^\pm_m(0)$ is equal to the coefficient of $q^2$ and is given by 
  $$
  \gamma_m^+=2 \sqrt{2}\sum_{j=0}^m \alpha_{2j}, \  \  \gamma_m^-=2 \sqrt{2}\sum_{j=0}^m \alpha_{2j+1}
 $$
 $(i)$~Let $n=2m$. Up to terms of order $>2$ one has
$$
\Delta_{n-1}^+/\Delta_{n}^+=\lambda^+_{m-1}(0)/\lambda^+_{m}(0)=\left(1+\gamma^+_{m-1}(q+q^2)\right)/\left(1+\gamma^+_{m}(q+q^2)\right)=
$$
$$
=1+(\gamma^+_{m-1}-\gamma^+_{m})q+(\gamma^+_{m-1}-\gamma^+_{m})(1-\gamma^+_{m}) q^2
$$
One has $\gamma^+_{m-1}-\gamma^+_{m}=-2 \sqrt{2}\,\alpha_{n}$, thus one obtains \eqref{qq3}. Similarly, one has, up to higher terms
$$
\Delta_{n+1}^-/\Delta_{n}^-=\lambda^-_{m}(0)/\lambda^-_{m-1}(0)=\left(1+\gamma^-_{m}(q+q^2)\right)/\left(1+\gamma^-_{m-1}(q+q^2)\right)=
$$
$$
=1+(\gamma^-_{m}-\gamma^-_{m-1})q+(\gamma^-_{m}-\gamma^-_{m-1})(1-\gamma^-_{m-1}) q^2
$$
One has $\gamma^-_{m}-\gamma^-_{m-1}=2 \sqrt{2}\,\alpha_{n+1}$, thus one obtains \eqref{qq3b}. 
 \newline
$(ii)$~Let $n=2m+1$. Up to terms of order $>2$ one has
$$
\Delta_{n+1}^+/\Delta_{n}^+=\lambda^+_{m+1}(0)/\lambda^+_{m}(0)=\left(1+\gamma^+_{m+1}(q+q^2)\right)/\left(1+\gamma^+_{m}(q+q^2)\right)=
$$
$$
=1+(\gamma^+_{m+1}-\gamma^+_{m})q+(\gamma^+_{m+1}-\gamma^+_{m})(1-\gamma^+_{m}) q^2
$$
One has $\gamma^+_{m+1}-\gamma^+_{m}=2 \sqrt{2}\,\alpha_{n+1}$, thus one obtains \eqref{qq4}. Similarly, one has, up to higher terms
$$
\Delta_{n-1}^-/\Delta_{n}^-=\lambda^-_{m-1}(0)/\lambda^-_{m}(0)=\left(1+\gamma^-_{m-1}(q+q^2)\right)/\left(1+\gamma^-_{m}(q+q^2)\right)=
$$
$$
=1+(\gamma^-_{m-1}-\gamma^-_{m})q+(\gamma^-_{m-1}-\gamma^-_{m})(1-\gamma^-_{m}) q^2
$$
One has $\gamma^-_{m-1}-\gamma^-_{m}=-2 \sqrt{2}\,\alpha_{n}$, thus one obtains \eqref{qq4b}. \endproof 
As a corollary we get the following expansions for the $a(n)^2$,
\begin{prop}\label{an2q2} $(i)$~Let $n=2m$. One has 
\begin{equation}\label{an2q21}	
	a(n)^2=(n+\frac 12)(n+1)\left(1+2 \sqrt{2}\,(\alpha_{n+1}-\alpha_{n})\ q+2 \sqrt{2}(\alpha_{n+1}-\alpha_{n})\ q^2+ b_nq^2+0(q^3)\right)
\end{equation}	
where $b_n=2 \sqrt{2}\,\alpha_{n}\gamma^+_{m}-2 \sqrt{2}\,\alpha_{n+1}\gamma^-_{m-1}-8\alpha_{n}\alpha_{n+1}$.\newline
$(ii)$~Let $n=2m+1$. One has 
\begin{equation}\label{an2q21}	
	a(n)^2=(n+\frac 12)(n+1)\left(1+2 \sqrt{2}\,(\alpha_{n+1}-\alpha_{n})\ q+2 \sqrt{2}(\alpha_{n+1}-\alpha_{n})\ q^2+ b_nq^2+0(q^3)\right)
\end{equation}	
where $b_n=2 \sqrt{2}\,\alpha_{n}\gamma^-_{m}-2 \sqrt{2}\,\alpha_{n+1}\gamma^+_{m}-8\alpha_{n}\alpha_{n+1}$.
\end{prop}
\proof $(i)$~By \eqref{nndet1} one has, for $n=2m$, 
$$
a(n)^2=(n+\frac 12)(n+1)\left(\Delta_{n-1}^+/\Delta_{n}^+\right)\left(\Delta_{n+1}^-/\Delta_{n}^-\right)
$$
and by Lemma \ref{consratio}, one gets that the product of the last terms gives
$$
1+2 \sqrt{2}\,(\alpha_{n+1}-\alpha_{n})\ q+2 \sqrt{2}\,\alpha_{n+1}(1-\gamma^-_{m-1})\ q^2-2 \sqrt{2}\,\alpha_{n}(1-\gamma^+_{m})\ q^2-8\alpha_{n}\alpha_{n+1}q^2
$$
which gives the value of $b_n$.\newline
$(ii)$~By \eqref{nndet1} one has, for $n=2m+1$,
$$
a(n)^2=(n+\frac 12)(n+1)\left(\Delta_{n+1}^+/\Delta_{n}^+\right)\left(\Delta_{n-1}^-/\Delta_{n}^-\right)
$$
and by Lemma \ref{consratio}, one gets that the product of the last terms gives
$$
1+2 \sqrt{2}\,(\alpha_{n+1}-\alpha_{n})\ q+2 \sqrt{2}\,\alpha_{n+1}(1-\gamma^+_{m})\ q^2-2 \sqrt{2}\,\alpha_{n}(1-\gamma^-_{m})\ q^2-8\alpha_{n}\alpha_{n+1}q^2
$$
which gives the value of $b_n$.\endproof 
We now need to determine the behavior of the expression 
$$
\beta_n:=(n+\frac 12)(n+1)b_n-n(n-\frac 12)b_{n-1}
$$
Let $n=2m+1$, then we get 
$$
(n+\frac 12)(n+1)b_n=(n+\frac 12)(n+1)\left(2 \sqrt{2}\,\alpha_{n}\gamma^-_{m}-2 \sqrt{2}\,\alpha_{n+1}\gamma^+_{m}-8\alpha_{n}\alpha_{n+1}\right)
$$
$$
n(n-\frac 12)b_{n-1}=n(n-\frac12)\left(2 \sqrt{2}\,\alpha_{n-1}\gamma^+_{m}-2 \sqrt{2}\,\alpha_{n}\gamma^-_{m-1}-8\alpha_{n-1}\alpha_{n}\right)
$$
We use \eqref{shiftan} to write all the terms in $\beta_n$ as multiples of $\alpha_n$. This gives as coefficient of $\alpha_n$ 
$$
2 \sqrt{2}\,(n+\frac 12)(n+1)\gamma^-_{m}+2 \sqrt{2}\,n(n-\frac 12)\,\gamma^-_{m-1}+2 \sqrt{2}\,(n+\frac 12)^2\,\gamma^+_{m}+2 \sqrt{2}\,n^2\,\gamma^+_{m}+
$$
$$
-8(n+\frac 12)(n+1)\alpha_{n+1}+8n(n-\frac 12)\alpha_{n-1}
$$
One has $\gamma^-_{m}=\gamma^-_{m-1}+2 \sqrt{2}\,\alpha_n$ and one can rewrite the coefficient as
\begin{equation}\label{betan}
\beta_n/\alpha_n=2 \sqrt{2}\,\left((2 n^2+n+\frac{1}{2})\,\gamma^-_{m}+(2 n^2+n+\frac{1}{4})\,\gamma^+_{m}\right)+8\,t_n
\end{equation}
$$
t_n=(n+\frac 12)(n+1)(-\alpha_{n+1})+n(n-\frac 12)(\alpha_{n-1}-\alpha_{n})
$$
Using \eqref{shiftan} one gets 
$$
8\,t_n=2 \left(-4 n^2+6 n+1\right)\,\alpha_n
$$
Let $n=2m$, then we get 
$$
(n+\frac 12)(n+1)b_n=(n+\frac 12)(n+1)\left(2 \sqrt{2}\,\alpha_{n}\gamma^+_{m}-2 \sqrt{2}\,\alpha_{n+1}\gamma^-_{m-1}-8\alpha_{n}\alpha_{n+1}\right)$$
$$
n(n-\frac 12)b_{n-1}=n(n-\frac12)\left(2 \sqrt{2}\,\alpha_{n-1}\gamma^-_{m-1}-2 \sqrt{2}\,\alpha_{n}\gamma^+_{m-1}-8\alpha_{n}\alpha_{n-1}\right)
$$
This gives as coefficient of $\alpha_n$ the expression
$$
2 \sqrt{2}\,(n+\frac 12)(n+1)\gamma^+_{m}+2 \sqrt{2}\,(n+\frac 12)^2\,\gamma^-_{m-1}
+2 \sqrt{2}\,n(n-\frac 12)\,\gamma^+_{m-1}+2 \sqrt{2}\,n^2\,\gamma^-_{m-1}+
$$
$$
-8(n+\frac 12)(n+1)\alpha_{n+1}+8n(n-\frac 12)\alpha_{n-1}
$$
One has $\gamma^+_{m}=\gamma^+_{m-1}+2 \sqrt{2}\,\alpha_n$ and one can rewrite the coefficient, with $t_n$ as above, as
\begin{equation}\label{betan1}
\beta_n/\alpha_n=2 \sqrt{2}\,\left((2 n^2+n+\frac{1}{2})\,\gamma^+_{m}+(2 n^2+n+\frac{1}{4})\,\gamma^-_{m-1}\right)+8\,t_n
\end{equation}
We thus obtain the expression of the first terms of the expansion of $-i\, d_n$,
\begin{lem}\label{idn} The first terms of the $q$-expansion of $-i\, d_n$ are
	\begin{equation}\label{midn}
-i\, d_n=-\frac{1}{\sqrt 2}(16 n^2+8 n +3)\alpha_n\, (q+q^2)+\alpha_n \, \rho_n\,q^2+O(q^3)
\end{equation}
where $\rho_n\in \Z[1/2]$ is given by
\begin{equation}\label{rhon} \rho_n=  (16 n^2+8 n+2)\sum_0^n \alpha_k+2 \left(-4 n^2+6 n+1\right)\alpha_n+2
\sum_{0\leq k\leq n, k=n (2)} \alpha_k
\end{equation}
\end{lem}
Now, as seen before Lemma \ref{detastrace}, one has
$$
\sum_0^n \alpha_k=\frac{(-1)^n \, _2F_1\left(1,n+\frac{3}{2};n+2;-1\right) \Gamma \left(n+\frac{3}{2}\right)}{\sqrt{\pi } \Gamma (n+2)}+\frac{1}{\sqrt{2}}
$$
which implies the equality 
\begin{equation}\label{hyperg}
\sum_0^{n-1} \alpha_k=- _2F_1\left(1,n+\frac{1}{2};n+1;-1\right) \alpha_n+\frac{1}{\sqrt{2}}
\end{equation}
This means that one can rewrite \eqref{rhon} as 
$$
\rho_n=\left(8 n^2+20 n+4-(16 n^2+8 n+2)_2F_1\left(1,n+\frac{1}{2};n+1;-1\right)\right)\alpha_n
$$
$$
+\frac{(16 n^2+8 n+2)}{\sqrt{2}}+2
\sum_{0\leq k\leq n, k=n (2)} \alpha_k
$$
In order to give a closed form to the last term, we use the Pochhammer symbol
$$(a)_k= \begin{cases}1 & \text { if } k=0 \\ a(a+1) \cdots(a+k-1) & \text { if } k>0\end{cases}
$$ 
and the equality 
\begin{equation}\label{signedsum}
\sum_0^n (-1)^k\alpha_k=\frac{3 \left(\frac{5}{2}\right)_{n-1}}{2 (2)_{n-1}}=(-1)^n(2n+1)\alpha_n
\end{equation}
which then gives, for $n$ even, using 
$$
2\sum_{0\leq k\leq n, k=n (2)} \alpha_k=\alpha_n+\sum_0^{n-1} \alpha_k+\sum_0^n (-1)^k\alpha_k
$$ 
the equality
\begin{equation}\label{hyperg1}
2\sum_{0\leq k\leq n, k=n (2)} \alpha_k=\left((2n+2)- _2F_1\left(1,n+\frac{1}{2};n+1;-1\right) \right)\alpha_n+\frac{1}{\sqrt{2}}
\end{equation}
This equality also holds for $n$ odd, using in that case
$$
2\sum_{0\leq k\leq n, k=n (2)} \alpha_k=\alpha_n+\sum_0^{n-1} \alpha_k-\sum_0^n (-1)^k\alpha_k
$$
Thus  one can rewrite the formula for $\rho_n$ as 
\begin{equation}\label{hyperg2}
\left(8 n^2+22 n+6-(16 n^2+8 n+3)_2F_1\left(1,n+\frac{1}{2};n+1;-1\right)\right)\alpha_n+\frac{(16 n^2+8 n+3)}{\sqrt{2}}
\end{equation}
We see that the last term cancells, after multiplication by $\alpha_n$, with the other term in $q^2$ in \eqref{midn} and we get the following result.
\begin{thm}\label{hinthyperg} The coefficient $\delta_2$ of $q^2$ in $-i\, d_n$ is given by
\begin{equation}\label{hyperg3}
\delta_2=\left(8 n^2+22 n+6-\left(16 n^2+8 n+3\right) \, _2F_1\left(1,n+\frac{1}{2};n+1;-1\right)\right)\alpha_n^2
\end{equation}	
One has for all $n$,
\begin{equation}\label{hyperg4}
-i\, d_n=-\frac{16 n^2+8 n+3}{\sqrt{2}}\,\alpha_n\,q+\end{equation} $$\left(8 n^2+22 n+6-\left(16 n^2+8 n+3\right) \, _2F_1\left(1,n+\frac{1}{2};n+1;-1\right)\right)\alpha_n^2\, q^2+O(q^3)
$$
\end{thm}
It is quite striking that in the formula for the term in $q^2$ there is no input of $\sqrt 2$ and such terms are provided by the values of the hypergeometric function whose first seven values are
$$
2-\sqrt{2},\,\ \frac{4}{3} \left(\sqrt{2}-1\right),\,-\frac{2}{5}  \left(4 \sqrt{2}-7\right),\,\frac{8}{35} \left(8 \sqrt{2}-9\right),\,-\frac{2}{63}  \left(64 \sqrt{2}-107\right)$$ $$\frac{4}{231} \left(128 \sqrt{2}-151\right),\,-\frac{2}{429}  \left(512 \sqrt{2}-835\right)
$$

\subsection{Coefficient of $q^3$ in $\det(\id +A_n)$	}
In order to compute the coefficient of $q^3$ in the expansion of $-i\,d_n$ we must first determine the coefficient of $q^3$ in the determinant of the matrices $\id +\cA_n^\pm$. It will be convenient to adopt the following notation 
\begin{equation}\label{An}
	A_n:=\begin{cases}
	\cA_n^+\ \textit{if}\ n=0 (2)\\
	\cA_n^- \     \textit{if}\ n=1 (2)	
	\end{cases}
\end{equation}

\begin{lem} \label{coeffq3} $(i)$~The $q$-expansion of  $\det(\id +A_n)$	is of the form $1+\sum c_j(n)q^j$ where
\begin{equation}\label{coeffq2j}
c_1(n)=c_2(n)=2 \sqrt{2}\sum_{0\leq j\leq n, j=n (2)} \alpha_{j}=\sqrt{2}\left((2n+2)- _2F_1\left(1,n+\frac{1}{2};n+1;-1\right) \right)\alpha_n+1
\end{equation}
$(ii)$~The coefficient $c_3(n)$ of $q^3$ is 
\begin{equation}\label{coeffq3j}
c_3(n)=2\sqrt 2\sum_{0\leq j\leq n, j=n (2)}\,\alpha_j\left(1-\frac 12 (4j+1)^2\right)
\end{equation}
$(iii)$~One has 
\begin{equation}\label{coeffq3jh}
c_3(n)=2 \sqrt{2}\left(-\frac 14 \, _2F_1\left(1,n+\frac{1}{2};n+1;-1\right)-3 ( n+\frac 12) \, _2F_1\left(1,n+\frac{3}{2};n+1;-1\right)\right)\alpha_n\end{equation} $$+\frac{ \sqrt{2}}{5} \left(-16 n^3-52 n^2-17 n+5\right)\alpha_n+2
$$
\end{lem}
\proof $(i)$~The determinant of $\id +A_n$ is expressed in terms of the sum of the traces on the exterior powers of $A_n$ and for the coefficient of $q^j$, $j\leq 3$, only the first wedge power enters 
as follows from Lemma \ref{keylem2}. The value of $c_1(n)=c_2(n)$ follows from Proposition \ref{coeffqeigen}, \eqref{gammanpm}, together with \eqref{hyperg1}.\newline
$(ii)$~We need to compute the term in $q^3$ in the $q$-expansion of the trace $\tr(\cA_n^\pm)$. By Proposition \ref{matrixentries}, and \eqref{gammaij}, the diagonal elements of the matrices  $\cA_n^\pm$ are given by 
$$
\gamma_{j,j}(q)=\sum_{\ell=0}^\infty (-1)^{j+\ell} 2^{\frac 32-2\ell-2j} \sigma(\ell,j)	\frac{q^{2\ell+1}}{1-q^{2\ell+1}}
$$
Only the two terms for $\ell=0,1$ contribute to the coefficient of $q^3$. One has 
by \eqref{sigmalm}
$$
	\sigma(\ell,m):=\binom{2\ell}{\ell}\binom{2m}{m}p_\ell(m)^2
	$$
while by \eqref{pmzg2}
$$
p_\ell(m):=\sum_{k=0}^\ell  2^{3 k}\frac{k! \prod _{j=0}^{k-1} \left(\ell-j\right)}{(2 k)! }\binom{m}{k}=\sum_{k=0}^\ell  2^{3 k}\binom{m}{k}\binom{\ell}{k}\binom{2k}{k}^{-1}
$$
so that $p_0(m)=1$, and $p_1(m)=4m+1$ for all $m$. Thus the first terms of the $q$-expansion of $\gamma_{j,j}(q)$ are 
$$
\gamma_{j,j}(q)=(-1)^{j}2^{\frac 32-2j}\binom{2j}{j}(q+q^2+q^3)-(-1)^{j}2^{\frac 32-2j-2}2\binom{2j}{j}(4j+1)^2q^3+ O(q^4)
$$
so that, using \eqref{alphal}, \ie $\alpha_\ell=(- 4)^{-\ell}\binom{2\ell}{\ell}$, the coefficient $c_3(n)$ of $q^3$ in the $q$-expansion of the trace $\tr(\cA_n^\pm)$, with the exponent $\pm $ given by the parity of $n$, is given by \eqref{coeffq3j}.\newline
$(ii)$~Now one computes the sums 
$$
\sum_0^n j\,\alpha_j=\frac{(-1)^n \left(n+\frac{1}{2}\right)! \, _2F_1\left(1,n+\frac{3}{2};n+1;-1\right)}{\sqrt{\pi } n!}-\frac{1}{4 \sqrt{2}}
$$
$$
\sum_0^n j^2\,\alpha_j=\frac{(-1)^n \, _2F_1\left(1,n+\frac{3}{2};n+1;-1\right) \Gamma \left(n+\frac{3}{2}\right)}{4 \sqrt{\pi } \Gamma (n+1)}+\frac{(-1)^n \Gamma \left(n+\frac{3}{2}\right)}{2 \sqrt{\pi } \Gamma (n)}-\frac{1}{16 \sqrt{2}}
$$
One has 
$$
\frac{(-1)^n \left(n+\frac{1}{2}\right)!}{\sqrt{\pi } \,n!}=\frac{(-1)^n \Gamma \left(n+\frac{3}{2}\right)}{\sqrt{\pi } \Gamma (n+1)}=(n+\frac 12) \alpha_n
$$
so the  above equalities can be rewritten as 
$$
\sum_0^n j\,\alpha_j=(n+\frac 12) \,_2F_1\left(1,n+\frac{3}{2};n+1;-1\right)\alpha_n-\frac{1}{4 \sqrt{2}}
$$
$$
\sum_0^n j^2\,\alpha_j=\frac 14 (n+\frac 12) \,_2F_1\left(1,n+\frac{3}{2};n+1;-1\right)\alpha_n+\frac 12 \,n(n+\frac 12) \alpha_n-\frac{1}{16 \sqrt{2}}
$$
Moreover one has 
$$
\sum_0^n (-1)^j j\,\alpha_j=(-1)^n \frac{n(2n+1)}{3}\alpha_n, \ \ 
\sum_0^n (-1)^j j^2\,\alpha_j=(-1)^n \frac{n(2n+1)(2n+3)}{15}\alpha_n
$$
Thus we get the equality, for $n$ even and $n$ odd 
$$
2\sum_{0\leq k\leq n, k=n (2)}k\, \alpha_k=\left(\frac{n(2n+1)}{3}+(n+\frac 12) \,_2F_1\left(1,n+\frac{3}{2};n+1;-1\right) \right)\alpha_n-\frac{1}{4 \sqrt{2}}
$$
and similarly 
$$
2\sum_{0\leq k\leq n, k=n (2)}k^2\, \alpha_k=\left(\frac{1}{60} n (2 n+1) (12 n+23)+\frac{(n+\frac 12)}{4} \,_2F_1\left(1,n+\frac{3}{2};n+1;-1\right) \right)\alpha_n-\frac{1}{16 \sqrt{2}}
$$
Combining these identities with \eqref{coeffq3j} gives \eqref{coeffq3jh}. \endproof 
\subsection{Jacobi coefficients $a(n)^2$ up to  $O(q^4)$}
In order to compute $a(n)^2$ up to terms in $O(q^4)$ we first prove a general fact concerning the expansion up to $O(q^4)$ of the product of ratios of the form $$R(n):=(a(n+1)/a(n-1))\times (a(n-2)/a(n))$$ where the $a(n)$ are given in terms of the variables $t(j)$ by 
$$
a(n):=1+s(n)(q+q^2), \ \ s(n):=\sum_{0\leq k\leq n, k=n (2)}t(k)
$$
\begin{lem} With the above notations, $R(n)$ is given, up to terms in $O(q^4)$, by
$$
1+q (t(n+1)-t(n))+q^2 (t(n) (s(n)-1)-t(n+1) (s(n-1)-1)-t(n+1) t(n))+$$ $$
q^3 \left(t(n+1) (s(n-1)-2) s(n-1)-t(n)(s(n)-2) s(n)+t(n) t(n+1) (s(n-1)+s(n)-2))\right)$$
\end{lem}
\proof This is checked by direct computation. \endproof

In our case one has $t(n)=2\sqrt 2\, \alpha_n$, and we can use \eqref{shiftan} in the form 
$$
(n-\frac 12)\,t(n-1)=-n\, t(n), \ \ (n+1)\,t(n+1)=-(n+\frac 12)\,t(n)
$$ 
This allows us to write, in our case, the coefficient of $q^3$ in the difference $$\eta_n:=(n+\frac 12) (n+1)R_n-(n-\frac 12) n\, R_{n-1}$$ as follows, with $z(n):=(s(n)-2) s(n)$ and $y(n):=s(n-1)+s(n)-2$,
$$
\eta_n=t(n)\left(-((n+\frac 12)^2+n^2)z(n-1)-(n+\frac 12) (n+1)z(n)-n(n-\frac 12)z(n-2)\right)+
$$
$$
+t(n)^2\left(n^2y(n-1)-(n+\frac 12)^2y(n)\right)
$$
\subsubsection{Terms in $z$.}
One has $s(n-2)=s(n)-t(n)$ and thus 
$$
z(n-2)=(s(n-2)-2) s(n-2)=(s(n)-t(n)-2)(s(n)-t(n))
$$
Moreover in our case we also have, by \eqref{signedsum}, $s(n)-s(n-1)=(2n+1)t(n)$, which gives
$$
z(n-1)=(s(n-1)-2) s(n-1)=(s(n)-(2n+1)t(n)-2)(s(n)-(2n+1)t(n))
$$
We can then write the contribution to $\eta_n$ of the terms in $z$ as 
$$
(-4 n^2-2 n-\frac{3}{4})\,t(n)\,s(n)^2+\frac{1}{2} \left(16 n^2+8 n+3+\left(16 n^3+20 n^2+4 n+1\right) t(n)\right)\,t(n)\,s(n)+
$$
$$
-\frac{1}{4} t(n)^2 \left(32 n^3+40 n^2+8 n+2+\left(32 n^4+48 n^3+32 n^2+6 n+1\right) t(n)\right)
$$
\subsubsection{Terms in $y$.}
We also get 
$$
y(n-1)=s(n-2)+s(n-1)-2=s(n)-t(n)+s(n)-(2n+1)t(n)-2=2(s(n)-(n+1)t(n)-1)
$$
and 
$$
y(n)=s(n-1)+s(n)-2=2s(n)-(2n+1)t(n)-2
$$
We can then write  $\eta_n$ as follows in terms of $t(n)$ and $s(n)$,
$$
\eta_n=-\frac 14(16 n^2+8 n+3)\,t(n)\,s(n)^2+\frac{1}{2} \left(4 (4 n+5) n^2 t(n)+16 n^2+8 n+3\right)\,t(n) \,s(n)
$$
\begin{equation}\label{dntotal}
-n^2 t(n)^2 \left(\left(8 n^2+12 n+7\right) t(n)+8 n+10\right)
\end{equation}
\subsubsection{Additional terms.}
They are given, using $p(n)=(1-\frac 12(4n+1)^2$, by 
$$
(n+\frac 12) (n+1)(t(n+1)p(n+1)-t(n)p(n))-(n-\frac 12) n\,(t(n)p(n)-t(n-1)p(n-1))
$$
which can be rewritten as 
$$
-(n+\frac 12)^2p(n+1)t(n)-(n+\frac 12) (n+1)t(n)p(n)-(n-\frac 12) n\,t(n)p(n)-n^2t(n)p(n-1)
$$
This gives the following additional contribution to the coefficient of $q^3$
\begin{equation}\label{additional}
\left(32 n^4+32 n^3+44 n^2+18 n+\frac{21}{8}\right)t(n)
\end{equation}
\subsection{Coefficient of $q^3$  in $-i\,d_n$}
We can now write the expression of the coefficient $\delta_3$ of $q^3$ in $-i\,d_n$ and it is more convenient to use the notation $t(n):=2\sqrt 2\, \alpha_n$, and to express the result in terms of the hypergeometric function $_2F_1\left(1,n+\frac{1}{2};n+1;-1\right)$ which was already involved in the coefficient of $q^2$ in Theorem \ref{hinthyperg} for the coefficient $\delta_2$ of $q^2$. 
\begin{thm}
	The coefficient $\delta_3$ of $q^3$ in $-i\,d_n(q)=\sum \delta_j q^j$ is given as a quadratic polynomial in terms of   $Y:=\, _2F_1\left(1,n+\frac{1}{2};n+1;-1\right)$ using the notation $t(n):=2\sqrt 2\, \alpha_n$, by the expression 
	$$
	\delta_3=-\frac{1}{16} \left(16 n^2+8 n+3\right) t(n)^3\,Y^2+\frac{1}{4} \left(4 n^2+11 n+3\right) t(n)^3\,Y$$ \begin{equation}\label{thmq3}-\left(4 n^4+4 n^3+\frac{23 n^2}{4}+\frac{7 n}{2}+\frac{3}{4}\right)t(n)^3+\left(32 n^4+32 n^3+48 n^2+20 n+\frac{27}{8}\right)t(n)
	\end{equation} 
\end{thm}
\proof This follows by adding to $\eta_n$ as in \eqref{dntotal}, the term \eqref{additional} and replacing $s(n)$ by its value given by \eqref{hyperg1}, \ie $s(n)=t(n) \left(n+1-\frac{Y}{2}\right)+1$.\endproof 
Note that the notation $t(n):=2\sqrt 2\, \alpha_n$ eliminates the $\sqrt 2$ in the formula \eqref{coeffqpert} which becomes 
\begin{equation}\label{thmq1}
\delta_1=-\frac{1}{4}(16 n^2+8 n +3)t(n)
\end{equation}
while the formula \eqref{hyperg3} for $\delta_2$ becomes
\begin{equation}\label{thmq2}
\delta_2=\frac 18 \left(8 n^2+22 n+6-\left(16 n^2+8 n+3\right) \,Y\right)t(n)^2
\end{equation}
Note also that the coefficient of the leading power of $Y$ in all these formulas  always involves the polynomial $\left(16 n^2+8 n+3\right)$ and that the coefficient of $Y$ in \eqref{thmq3} involves $\frac{1}{4}\left(4 n^2+11 n+3\right)$ which also appears as $\frac 18\left(8 n^2+22 n+6\right)$ as the coefficient of $Y^0$ in \eqref{thmq2}.


\end{document}